\documentclass[lettersize,journal]{IEEEtran}
\usepackage{amsmath,amsfonts}
\usepackage{algorithmic}
\usepackage{algorithm}
\usepackage{array}
\usepackage[caption=false,font=normalsize,labelfont=sf,textfont=sf]{subfig}
\usepackage{textcomp}
\usepackage{stfloats}
\usepackage{url}
\usepackage{xr} 
\usepackage{verbatim}
\usepackage{graphicx}
\usepackage{cite}
\usepackage{longtable} 
\usepackage{makecell}
\usepackage{hyperref}
\newcommand{\R}{\mathbb R}
\newtheorem{theorem}{Theorem}[section]
\newtheorem{proposition}{Proposition}[section]% 
\newtheorem{lemma} {Lemma}[section]
\newtheorem{remark} {Remark}[section]%
\newtheorem{definition} {Definition}[section]%
\newtheorem{assumption}{Assumption}[section]
\begin{document}

\title{Newton Method for Soft Quadratic Surface Support Vector Machine with 0-1 Loss Function}

\author{Guoping Li, Wen Song*
	
\thanks{School of Mathematical Sciences, Harbin Normal University, Harbin, 150025,  P.R. China email:18845562781@163.com; wsong@hrbnu.edu.cn}
\thanks{* Corresponding author}}

        % <-this % stops a space



\maketitle

\begin{abstract}
A nonlinear kernel-free soft quadratic surface support vector machine model with 0-1 loss function ($L_{0/1}$-SQSSVM) is proposed  for  binary classification problems, which is non-convex discontinuous.  We are devoted to establishing the first  and the second-order optimality conditions for the $L_{0/1}$-SQSSVM.  We  establish a stationary equation using the properties of proximal operator of 0-1 loss function. We design a  Newton method based on the stationary equation to solve $L_{0/1}$-SQSSVM model and prove that the Newton method has local quadratic convergence under the second-order sufficient condition. Numerical experience on artificial datasets and benchmark datasets demonstrate that the Newton method for $L_{0/1}$-SQSSVM achieves higher classification accuracy with less CPU time cost than other state-of-the-art	methods.
\end{abstract}

\begin{IEEEkeywords}
0-1 loss function, soft quadratic surface support vector machine, subderivatives, optimality condition,  Newton method.
\end{IEEEkeywords}

\section{Introduction}
\IEEEPARstart{S}{upport} vector machines (SVM) were first introduced by Cortes and Vapnik \cite{CV} and have been widely adopted for pattern recognition \cite{CJB}, data analysis and classification \cite{CNSP, APR},  financial prediction \cite{Cao}, texture classification \cite{sharif2019optimised, Kim}, and so on. In binary classification problems, the main idea of SVM is to find a hyperplane that separates the training datasets into two classes.  Given a linearly separable training data set $\{(x^1, y_1), \cdots, (x^i, y_i), \cdots, (x^n, y_n)\}$, where $x^i \in \R^m$ indicates the position of the $i$-th training point and $y_i \in \{-1, +1\}$ is its label. The objective is to find a hyperplane $h(x) = \langle w,x\rangle +b$ that  maximizes the margin to separate the sample points  as distinctly as possible. The parameters of  this hyperplane are obtained by solving the following optimization problem \cite{Boser}:
\begin{align}\label{SVM}
	\underset{w\in \mathbb{R}^{m},b \in \R}{\rm min}~&\frac{1}{2}\|w\|^2_2 \notag\\
	{\rm s.t.~~} ~& y_i(\langle w,x^i\rangle +b)\ge 1,~i=1,2,\cdots,n.
\end{align}
This model is known as the hard-margin SVM, as it requires the correct classification of all training samples. However, most of training datasets are not completely separated by linear hyperplane in real applications. The popular method is to allow violations in the satisfaction of the constraints in \eqref{SVM} and penalize such violations in the objective function, the model is called the soft margin SVM model (SSVM):
\begin{align*}
	\underset{w\in \mathbb{R}^{m},b \in \R}{\rm min}~&\frac{1}{2}\|w\|^2_2 + \lambda \underset{i=1}{\overset{n}{\sum}} \ell(1-y_i(\langle w,x^i\rangle +b)), 
\end{align*} 
where $\lambda>0$ is a penalty parameter and $\ell(\cdot)$ is a  non-negative real-valued loss function. The loss function  is an important research topic in SSVM \cite{HMN}, some popular loss functions include the pinball loss \cite{Huang, tanveer2021sparse}, the hinge soft-margin loss \cite{Hinge}, the ramp soft-  margin loss \cite{ Ramp, carrizosa2014heuristic}, the least squares loss \cite{suykens1999least, gao2025kernel}, and so on. Cortes and Vapnik \cite{CV}, Brooks \cite{Brooks}, Feng, Yang and Huang \cite{Feng} have pointed out that the ideal SSVM is
\begin{align*}
	\underset{w\in \mathbb{R}^{m},b \in \R}{\rm min}~&\frac{1}{2}\|w\|^2_2 + \lambda \underset{i=1}{\overset{n}{\sum}} \ell_{0/1}(1-y_i(\langle w,x^i\rangle +b)),
\end{align*} 
where $\ell_{0/1}(t) = \left\{\begin{array}{clc}
	1, &t>0,\\
	0, &t\leq 0
\end{array}
\right. $ is 0-1 loss function.  For the given training datasets, it is not always to be linearly separated. To address nonlinear separation, Cortes and Vapnik \cite{CV} $(\xi_i =0, i=1\ldots,n)$  and Vapnik \cite{Vapnik1998} then proposed a version with a nonlinear kernel function  $\phi : \R^m \to \R^l$ : 
\begin{align*}
	\underset{v\in \mathbb{R}^{l},b \in \R}{\rm min}~&\frac{1}{2}\|v\|^2_2 + \hat{\eta}\underset{i=1}{\overset{n}{\sum}} \xi_i\notag\\
	{\rm s.t.~~} ~& y_i(\langle v,\phi(x^i)\rangle +b)\ge 1- \xi_i,~i=1,2,\cdots,n.
\end{align*}
The idea of the kernel-based SVM model is to map data points to a higher dimensional space using  a nonlinear kernel function and then apply the SVM model for classification in this feature space. The well-known kernel functions include the  Gaussian kernel  and the quadratic kernel \cite{BSAJS}. Liu et al. \cite{liu2024nonlinear} skillfully integrate kernel methods with the 0-1 loss function to propose a nonlinear $L_{0/1}$-SVM model. However, There is no general principle for selecting an appropriate kernel function. Moreover, the performance of kernel-based SVM models heavily depends on the kernel parameters \cite{BSAJS, NCJS} and tuning those parameters often requires significant effort. 

To avoid the complexity and computational cost associated with the nonlinear kernel trick, a kernel-free quadratic surface support vector machines (QSSVM) model was proposed by Dagher \cite{DI}.  This model aims to separate data using a quadratic surface $h(x) = \frac{1}{2}x^TWx + b^Tx+c,$ where $ W\in \mathbb{S}^{m}$ is an $m\times m$ symmetric matrix,  $ b \in \R^m$ and $c \in \R$. The QSSVM optimization problem is formulated as follows: 
\begin{align}\label{QSSVM}
	\underset{W,b,c}{\rm min}~&\underset{i=1}{\overset{n}{\sum}}\frac{1}{2}\|Wx^i+b\|_2^2\notag\\
	{\rm s.t.} ~& y_i(\frac{1}{2}{x^i}^TWx^i+b^Tx^i+c)\ge 1,~i=1,2,\cdots,n,\\
	~& W=W^T\in \mathbb{R}^{m\times m}, ~x^i, b \in \R^m, ~c \in \R.\notag
\end{align}
Luo et al. \cite{Luo} developed an extension that incorporates noise and outliers, called the soft margin quadratic surface support vector machines (SQSSVM):
\begin{align*}
	\underset{W,b,c}{\rm min}~&\underset{i=1}{\overset{n}{\sum}}\frac{1}{2}\|Wx^i+b\|^2_2 + \hat{\eta} \underset{i=1}{\overset{n}{\sum}} \xi_i \notag\\
	{\rm s.t.} ~& y_i(\frac{1}{2}{x^i}^TWx^i+b^Tx^i+c)\ge 1 - \xi_i,~i=1,2,\cdots,n,
\end{align*}
where $\hat{\eta}$ is a penalty parameter and $\xi_i$ is a slack variable for each constraint. Following this, some scholars have studied the SQSSVM model with different loss functions. Bai et al. \cite{Bai} proposed a quadratic kernel-free least squares support vector machines for classification problems  and regression problems \cite{Ye, He2024quadratic}:
\begin{align*}
	\underset{W,b,c, e}{\rm min}~&\underset{i=1}{\overset{n}{\sum}}\|Wx^i+b\|^2_2 + \frac{C}{2} \underset{i=1}{\overset{n}{\sum}} e_i^2 \notag\\
	{\rm s.t.} ~& \frac{1}{2}{x^i}^TWx^i+b^Tx^i+c = y_i- e_i,~i=1,2,\cdots,n.
\end{align*}
This model penalizes the residual using a quadratic function and transforms inequality constraints into an equality constraint.  Mousavi et al. \cite{Mousavi2022}  propose a QSSVM with  $\ell_1$-norm regularization (L1-QSSVM):
\begin{align*}
	\underset{W,b,c}{\rm min}~&\underset{i=1}{\overset{n}{\sum}}\frac{1}{2}\|Wx^i+b\|^2_2 + \lambda \underset{1\le i\le j\le m}{\sum} |W_{ij}| \notag\\
	{\rm s.t.} ~& y_i(\frac{1}{2}{x^i}^TWx^i+b^Tx^i+c)\ge 1 ,~i=1,2,\cdots,n,
\end{align*}
where $\lambda>0$ and penalizes nonzero elements of the matrix $W$. Recently, Wu, Yang and Ye \cite {Yang} propose  a kernel-free  quadratic surface support vector machine model with 0-1 loss (QSSVM$_{0/1}$): 
\begin{align}\label{QSSVM0/1}
	\underset{W,b,c}{\rm min}~\underset{i=1}{\overset{n}{\sum}}\frac{1}{2}\|Wx^i+b\|^2+C\underset{i=1}{\overset{n}{\sum}}\ell_{0/1}(1-y_i(\frac{1}{2}{x^i}^TWx^i+b^Tx^i+c)).
\end{align}
Instead of considering the matrix version of QSSVM, the literatures \cite{Luo}-\cite{Yang} consider corresponding equivalent vector form. For example, the authors in \cite{Yang} derive the following equivalent  vector form of optimization  problem  \eqref{QSSVM0/1} 
\begin{align}\label{vecQSSVM0/1}
	\underset{\tilde{w}, b, c,u}{\rm min}~&\underset{i=1}{\overset{n}{\sum}}\frac{1}{2}\|M_i\tilde{w}+b\|^2+C\|u_+\|_0& \notag\\
	{\rm s.t.} ~& u+A\tilde{w}+ Bb +cy =1,~i=1,2,\cdots,n,\\
	~&\tilde{w} \in\R^{ \frac{m^2+m}{2}}, ~b \in \R^m,~c \in \R~,~ u \in \R^n,\notag
\end{align}
where $M_i \in\R^{m \times \frac{m^2+m}{2}}$, $A \in  \R^{n \times \frac{m^2+m}{2}}$, $B \in \R^{n \times m}$ and $\|u_+\|_0 = \underset{i=1}{\overset{n}{\sum}} \ell_{0/1}(u_i)$. \cite{Luo}-\cite{Mousavi2022} have proved the existence and uniqueness of the optimal solution and  designed optimization algorithms based on the vector form. \cite{Yang} proposed an alternating direction method of multipliers (ADMM) for problem \eqref{vecQSSVM0/1} using the proximal operator of 0-1 loss function.

The optimization problem with 0-1 loss function is NP-hard \cite{BK, EAVK} since the 0-1 loss function is nonconvex and discontinuous. It is challenging to design a efficient numerical algorithm for this problems since traditional optimization methods are not equipped to handle this type of problem. Relevant researches can be divided into two classes: the first class comprises methods based on continuous surrogates of the $\ell_{0}$-norm, such as $\ell_{1}$-norm  penalty, smoothing clipped absolute deviation penalty and capped $\ell_1$ penalty. The second class involves methods that directly optimize the 0-1 loss function \cite{WSZZX, Yang, ZPXQ, zhang2024zero}. As mentioned in \cite{WSZZX, Yang},   first-order methods like the ADMM and the inexact augmented Lagrangian method \cite{zhang2024zero} were proposed based on P-stationary point defined in terms of the proximal operator of the 0-1 loss function. Recently, some second-order methods have been designed to solve problem \eqref{vecQSSVM0/1}, such as the prime-dual active set method \cite{PD} and the Newton method \cite{{ZPXQ}, zhang2023inalm}. It is well known that the second order sufficient condition (positive definite of Hessian) plays the key role in deriving local quadratic convergence of Newton methods. However, there is no discussions about the second-order optimilaty conditions in the mentioned literatures. This motivated us in this paper to explore the second-order optimization conditions and based on them to derive local quadratic convergence of Newton methods.

In this paper, we consider the following general sparsity optimization problem:
\begin{align}\label{Prim}
	\underset{X\in \mathbb{R}^{m\times p}}{\rm min}~&f(X) + \lambda\underset{i=1}{\overset{n}{\sum}} \ell_{0/1}(1+\langle A_i, X\rangle),
\end{align}
where $f:\mathbb{R}^{m\times p}\to \mathbb{R}$ is twice continuously differentiable function and $\lambda > 0$ is a penalty parameter, $A_i \in \mathbb{R}^{m\times p}$ and $\langle  A_i, X \rangle = {\rm tr}(A_i^TX)$ is the inner product of matrix. If $p=1$, model  \eqref{Prim} reduces to the optimization problem in \cite{ZPXQ, ZXL}. If the matrices $X$ and $A_i$ take the following forms: 
$$X = \begin{pmatrix}
	W, {\rm Diag}(b),ce_1
\end{pmatrix}\in \R^{m \times (2m+1)}$$
and
$$A_i = -y_i\begin{pmatrix}
	\frac{1}{2}{x^i}{x^i}^T,
	{\rm Diag}(x^i),
	e_1
\end{pmatrix}\in \R^{m\times(2m+1)},$$ 
where ${\rm Diag}(b) \in\R^{m \times m} $ is a diagonal matrix with diagonal elements equal to $b$,  $e_1 = (1,0, \cdots, 0)^T\in \R^m$, we have 
$$1+\langle A_i, X\rangle = 1-y_i(\frac{1}{2}{x^i}^TWx^i+b^Tx^i+c), ~ i = 1,2, \cdots, n.$$ 
Hence problem \eqref{Prim}  is reduce to  \eqref{QSSVM0/1} if $f(X) =\underset{i=1}{\overset{n}{\sum}}\frac{1}{2}\|X\tilde{x}^i\|^2,$ where $\tilde{x}^i = \begin{pmatrix}
	x^i\\
	1_m\\
	0
\end{pmatrix} \in \R^{2m+1}$ and $1_m =(1, \cdots, 1)\in \R^m.$ The distinction is that the variable to be optimized in \eqref{vecQSSVM0/1} is  vectors, whereas it is a matrix in problem \eqref{Prim}. The main contributions of this paper  are as follows:\\
(a) For the composite  optimization problem \eqref{Prim}, we  derive the first and second-order optimality conditions using the first and second-order subderivatives. Moreover, we establish the relationship between the P-stationary points and the minimizer of problem \eqref{Prim} by using the properties of the proximal operator of 0-1 loss function. \\
(b) We propose a  Newton method for problem \eqref{Prim} based on the stationary equation. Moreover, we establish its local quadratic convergence under the second-order sufficient condition.\\
(c)  We test the performance of our method by conducting various numerical experiments on artificially generated data with two features and public benchmark data sets with more than two features. Extensive experimental results show that our method is more efficient than other state-of-the-art methods for binary classification problems.

	\section{Preliminaries}
Throughout the paper, for a given set $I \subset \{1,2,\ldots, n \}$, $|I|$ denotes the cardinality of  $I$ and $\bar{I}$ denotes its complement. For a vector $ z \in \mathbb{R}^n$,  $z_I$  represents the sub-vector  containing the elements of $z$ indexed on $I$.  $\|z_+\|_0$ counts the number of the positive  elements of $z$. $\|z\|$ denotes the Euclid norm of $z$. For a matrix $X\in \mathbb{R}^{m\times p}$, $\sigma(X)$ indicates the singular values vector of $X$. $\|X\|_F$ denotes  the Frobenius norm of $X$ and $\|X\|$ denotes  the spectral norm of $X$. The open neighborhood of $X^*$ with a radius $\delta>0$ is denoted by 
$B(X^*, \delta) = \{X\in \mathbb{R}^{m\times p}~|~\|X- X^*\|_F<\delta\}.$
The  proximal operator \cite[Definition 1.22]{RW} of a function $g(\cdot) $  with a constant $\alpha > 0$ is defined by
$${\rm prox}_{\alpha g} (z) := \underset{u }{\rm argmin} ~g(u) + \frac{1}{2\alpha}\|u-z\|^2.$$

\begin{definition} \cite{RW}
	Given a function $f \colon \mathbb{R}^n\to \bar{\mathbb{R}}$ and $\bar x  \in {\rm dom}f$  with $f(\bar{x})$ finite, the regular subdifferential of $f$ at $\bar{x}$ is defined by 
	\begin{equation*}
		\hat{\partial}f(\bar{x})=\{v \in \mathbb{R}^n~|~\liminf_{x\to \bar{x}}\frac{f(x)-f(\bar{x})-\langle v,x-\bar{x}\rangle }{\|x-\bar{x}\|}\geq 0\}.
	\end{equation*}
	The limit subdifferential of $f$ at $\bar{x}$ is defined by 
	\begin{equation*}
		\partial f(\bar{x})=\{v \in \mathbb{R}^n~|~\exists~x_k \overset{f}{\to}\bar{x},~v_k \in \hat{\partial}f(x_k),~{\rm s.t.}~v_k \to v\},
	\end{equation*}
	where $x_k \overset{f}{\to}\bar{x}$ stands for $x_k \to\bar{x}$ and $f(x_k)\to f(\bar{x})$.\\
	The horizon subdifferential of $f$ at $\bar{x}$ is defined by 
	\small{\begin{equation*}
		\partial^{\infty} f(\bar{x})=\{v \in \mathbb{R}^n|\exists~x_k \overset{f}{\to}\bar{x},~v_k \in \hat{\partial}f(x_k),~\lambda_k\downarrow 0,~{\rm s.t.}~\lambda_kv_k \to v\}.
	\end{equation*}}
\end{definition}

\begin{definition} \cite{Twice}\label{prox-subdiff}
	We say that $v \in \mathbb{R}^n $ is a proximal subgradient of $f$ at $\bar x$ if there exist $r \in \mathbb{R}_+$ and a neighborhood U of $\bar x$ such that for all $x \in  U$, we have
	$$ f(x)\ge f(\bar x) + \langle v, x-\bar x\rangle-\frac{r}{2}\|x-\bar x \|^2.$$  
	The set of all such $v$ is called the proximal subdifferential of $f$ at $\bar x$ and is denoted by $\partial^p f(\bar{x})$.
\end{definition}

\begin{definition} \cite{RW}
	For a function $f \colon \mathbb{R}^n\to \bar{\mathbb{R}}$ and a point $\bar{x}$ with $f(\bar{x})$ finite and $\partial f(\bar{x}) \neq \emptyset$, one has $f$ regular at $\bar{x}$ if and only if $f$ is locally lsc at $\bar{x}$ with 
	$$\partial f(\bar{x})= \hat{\partial} f(\bar{x}),\ \ \partial^{\infty} f(\bar{x}) = \hat{\partial} f(\bar{x})^{\infty},$$
	where $\hat{\partial} f(\bar{x})^{\infty} = \{v~|~\exists~ v^k \in \hat{\partial} f(\bar{x}),~\lambda^k\downarrow 0,~ {\rm with}~ \lambda^kv^k \to v\}$ is the horizon cones of $\hat{\partial} f(\bar{x})$.
\end{definition}

\begin{definition} \cite{RW}
	Given a function $f \colon \mathbb{R}^n\to \bar{\mathbb{R}}$ and a point $\bar{x} \in \mathbb{R}^n$ with $f(\bar{x})$ finite, the subderivative of $f$ at $\bar x$ is defined by 
	\begin{equation*}
		{\rm d} f(\bar{x})(w)= \liminf_{\substack{t \downarrow 0 \\ w^{'}\to w}}\frac{f(\bar{x}+t w^{'})-f(\bar{x})}{t},~w\in \mathbb{R}^n.
	\end{equation*}
\end{definition}

We observe by \cite[Exercise 8.4]{RW} that the regular subdifferential has a dual representation 
\begin{equation}\label{diff-der}
	\hat{\partial}f(\bar{x})=\{v \in \mathbb{R}^n~|~\langle v, w \rangle \le {\rm d}f(\bar x)(w) \mbox{ for all } w \in \mathbb{R}^n\}.
\end{equation}	
For any $\bar{x} \in {\rm dom}f$ and $\bar{v}\in \partial f(\bar x)$, the critical cone of $f$ at $\bar{x}$ is defined by 
\begin{equation}\label{critical}
	K_f(\bar{x},\bar{v}):=\{w \in\mathbb{R}^n~|~\langle \bar{v},w\rangle={\rm d} f(\bar{x})(w)\}.
\end{equation}

For a family of functions $f_t: \mathbb{R}^n\to \mathbb{R}$, we say that $f_t$ epi-converge to a function $f$ as $t\downarrow 0$, written as $f={\rm epi}${-}${\rm lim}_{t\downarrow 0}f_t$, if and only if the sets ${\rm epi}f_t$ converge to ${\rm epi}f$ as $t\downarrow 0$ in the Painlev\'{e}-Kuratowski sense.
For some $\bar x \in {\rm dom} f$ and $h,w\in \mathbb{R}^n$. Define two forms of second-order difference quotients for $f$ at $\bar x $ for $h$.
\begin{equation*}
	\Delta_t^2 f(\bar x)(h,w) := \frac{2}{t^2}(f(\bar{x}+th+\frac{1}{2}t^2w)-f(\bar x)-t{\rm d}f(\bar{x})(h))
\end{equation*}
and
\begin{equation*}
	\Delta_t^2 f(\bar x,v)(h) := \frac{2}{t^2}(f(\bar x+th)-f(\bar x)-t\langle v,h\rangle),
\end{equation*}
where $t>0$.

\begin{definition}\cite{RW}
	Assume that $f(\bar{x})$ and the subderivative ${\rm d}f(\bar x)(h)$ are finite, the parabolic subderivative of $f$ at $\bar x$ for $h$ in the direction $w$ is defined by 
	\begin{equation}
		{\rm d}^2f(\bar x)(h,w):= \liminf_{{t\downarrow 0}\atop {w^{'} \to w}}  \Delta_t^2 f(\bar x)(h,w^{'}).
	\end{equation}
\end{definition}

If the functions $\Delta_t^2 f(\bar x)(h,\cdot)$ epi-converge to ${\rm d}^2f(\bar x)(h,\cdot)$ as $t\downarrow 0$ and ${\rm d}^2f(\bar x)(h,\cdot) \not\equiv \infty$,  we say that $f$ is parabolically epi-differentiable at $\bar x$ for $h$.
\begin{definition} \cite{RW}
	For a function $f \colon \mathbb{R}^n\to \bar{\mathbb{R}}$, any $\bar x\in \mathbb{R}^n$ with $f(\bar x)$ finite, the second subderivative of $f$ at $\bar x$ for $v$ in a direction $h$ is defined as 
	\begin{equation*}
		{\rm d}^2f(\bar x,v)(h):= \liminf_{\substack{t \downarrow 0 \\ h^{'}\to h}} \Delta_t^2 f(\bar x,v)(h^{'}).
	\end{equation*}
\end{definition}

If the functions $\Delta_t^2 f(\bar x,v)(\cdot)$ epi-converge to ${\rm d}^2f(\bar x,v)(\cdot)$ at $\bar{x}$ as $t\downarrow 0$, we say that $f$ is twice epi-differentiable at $\bar{x}$ relative to $v\in \mathbb{R}^n$. If the function $f$ is twice epi-differentiable at $\bar{x}$ relative to $v$, then the second subderivative ${\rm d}^2f(\bar x,v)(\cdot)$ is lower semicontinuous and positively homogeneous of degree $2$. 
%\begin{lemma}
%For	$f:\R^n \to \bar{\R}$, any point $\bar{x}$ with $f(\bar{x})$ is finite and any vector $h$ %with ${\rm d}f(\bar{x})(h)$ finite, the function ${\rm d}^2f(\bar{x})(h, \cdot)$  is lsc. For %vectors $v$ with $\langle v,h\rangle={\rm d} f(\bar{x})(h)$ it satisfies 
%\end{lemma}
From \eqref {critical}, we have that  for any $\bar x \in {\rm dom}f$ and $v \in \partial f(\bar x)$, for every $h \in K_f(\bar x,v)$, we have ${\rm d} f(\bar{x})(h) = \langle v,h\rangle$. As shown in \cite[Proposition 13.64]{RW},  we have 
\begin{equation} \label{para-inequa}
	{\rm d}^2f(\bar x,v)(h)  \le \inf_{w \in \mathbb{R}^n} \{ {\rm d}^2f(\bar x)(h,w)- \langle v, w \rangle \}
\end{equation}
holds.

\begin{definition}\cite{RW}
	A function $f:\R^n \to \bar{\R}$ is parabolically regular at a point $\bar{x}$ for a vector $v$ if $f(\bar{x})$ is finite and the inequality in \eqref{para-inequa} holds with equality for every $h$ having ${\rm d}f(\bar{x})(h)=\langle v,h\rangle$, i.e.,
	\begin{equation} \label{para-equa}
		{\rm d}^2f(\bar x,v)(h)  =\inf_{w \in \mathbb{R}^n} \{ {\rm d}^2f(\bar x)(h,w)- \langle v, w \rangle \}.
	\end{equation}
\end{definition}

\begin{lemma} \cite{L}\label{lemma2.2}
	Let $g(y)= \|y_+\|_0 = \underset{i=1}{\overset{n}{\sum}}\ell_{0/1}(y_i)$, for given $ y^*\in \mathbb{R}^n$, the function $g$ is regular at $y^*$ and 
	\begin{align}\label{regular}
		\partial^pg( y^*)&=\hat{\partial}g( y^*)=\partial g( y^*)=\partial^{\infty}g( y^*) \nonumber\\
		&=\{v \in  \mathbb{R}^n ~|~v_i \ge 0,~{\rm if }~y_i^*=0;~v_i =0,~{\rm if }~y^*_i \neq 0)\}.
	\end{align}
	Moreover,
	\begin{equation*}
		{\rm d}g( y^*)(v)=0, ~v \in {\rm dom\,d}g(y^*),
	\end{equation*}
	where 
	\begin{equation*}
		{\rm dom~d}g( y^*) = \{v \in \mathbb{R}^n~|~ v_i \le 0,~{\rm if }~y_i^*=0,~v_i \in \mathbb{R},~{\rm if }~ y^*_i \neq 0 \}.
	\end{equation*} 	
\end{lemma}

\begin{lemma}\cite{L} \label{para subderivative} 
	The function $g$ is parabolically epi-differentiable at $ y^*$ for  $h \in {\rm dom~d}g(y^*)$ in the direction $w$ and
	\begin{align}\label{para-diff}
		{\rm d}^2g( y^*)(h,w) =0, ~~w \in {\rm dom\,d^2}g( y^*)(h,\cdot),
	\end{align}
	where 
	\begin{equation*}
		{\rm dom\,d^2}g(y^*)(h,\cdot) = \left\{w \in \mathbb{R}^n~\bigg|\begin{aligned}
		&w_i \le 0, ~{\rm if}~ y^*_i = 0, ~h_i = 0 ;\\
		&w_i \in\mathbb{R},~ {\rm if}~ y^*_i = 0, ~h_i \neq 0 
		\end{aligned}\right\}.
	\end{equation*}
\end{lemma}

\begin{lemma}\cite{L} \label{Prop 3.4}
	For given $y^* \in \mathbb{R}^n$ and $ v\in \partial g( y^*)$, the function $g$ is parabolically regular at $y^*$ for $v$,  and  the infimum in \eqref{para-equa} can be attained at some $ w\in {\rm dom\,d^2}g(y^*)(h,\cdot)$. Moreover, 
	\begin{align*}
		{\rm d}^2g(y^*,v)(h) =\left\{
		\begin{array}{clc}
			0,   & h \in K_g(y^*,v),\\
			+\infty, & {\rm otherewise},
		\end{array}
		\right.
	\end{align*}
	where 
	\begin{equation*}
		K_g( y^*,v) = \left\{u \in \mathbb{R}^n~\bigg|~
		\begin{aligned}
			~&u_i \in\mathbb{R},~y^*_i \neq 0 ;\\
			~&u_i \le 0,~y^*_i = 0, v_i =0 ;\\
			~&u_i =0,~y^*_i = 0, v_i \neq 0
		\end{aligned}
		\right\}.
	\end{equation*}
\end{lemma}

\section{The optimality conditions}\noindent
In this section, we  firstly calculate the first and second-order subderivatives of the  composite function. Then we establish the optimality conditions of problem \eqref{Prim} by using the relationship between subderivative and minimizer in variational analysis \cite{RW}.

\subsection{The first and second-order subderivatives of the composite function}\label{subsec1}

In this part, we  calculate the first and  second-order subderivatives of the  composite function by using some tools in convex analysis \cite{Fundatnentals}, variational analysis \cite{RW, Mordu} and perturbation analysis \cite{BS}.
We define  the mapping $F:\mathbb{R}^{m\times p}\to \mathbb{R}^{n}$ by
$$F(X)  = (F_1(X), F_2(X), \cdots, F_n(X))^T \in \R^n,$$ 
where $F_i(X) = 1+\langle A_i, X\rangle, ~i =  1,2, \cdots ,n.$
Then $F(X)$ is an affine mapping.  For convenience, denoted by $\varphi(X) = \|F(X)_+\|_0 = \underset{i=1}{\overset{n}{\sum}}\ell_{0/1}(1+\langle A_i, X\rangle) $. Let  
$$A = \begin{pmatrix}
	A_1\\
	A_2\\
	\vdots\\
	A_n
\end{pmatrix} \in \R^{nm\times p}$$ 
be a block matrix. Let $ \mathcal{T}_*:= \{i \in [n]~|~ F_i(X^*) = 0\}$ and $A_{\mathcal{T}_*}$ be the sub-matrix  containing the block of  $A$  indexed on $\mathcal{T}_*$. We define two types of block matrix multiplication, for a vector $v \in \R^n$ and a matrix $V \in \R^{m\times p}$,
\begin{align*}
	A^Tv &= \underset{i=1}{\overset{n}{\sum}} v_i A_i \in \R^{m\times p},\\
	AV &= \begin{pmatrix}
		\langle A_1, V \rangle\\
		\langle A_2, V\rangle\\
		\vdots\\
		\langle A_n, V\rangle
	\end{pmatrix} \in \R^{n}.
\end{align*}
For a matrix $X =(x_{ij}) \in \mathbb{R}^{m\times p}$, ${\rm vec}(X) = (x_{11}, \cdots, x_{m1}, \cdots, x_{1p}, \cdots, x_{mp})^T \in \R^{mp}$ represents the matrix vectorization. Then we have 
$\|X\|_F = \sqrt{\underset{i=1}{\overset{m}{\sum}}\underset{j=1}{\overset{p}{\sum}}|x_{ij}|^2} = \|{\rm vec}(X)\|$ and  for any $X_1, X_2 \in  \mathbb{R}^{m\times p} $,
\begin{equation}\label{FnormSubtraction}
	\|X_1-X_2\|_F = \|{\rm vec}(X_1)-{\rm vec}(X_2)\| .
\end{equation}
It follows from \cite[(2.3.7)]{golub2013matrix} and \eqref{FnormSubtraction}, we have that
\begin{equation}\label{spectral-norm}
	\|X_1-X_2\|\le \|X_1-X_2\|_F = \|{\rm vec}(X_1)-{\rm vec}(X_2)\|.
\end{equation}
Refer to \cite{HM}, the Hessian matrix of a differentiable mapping $f(X): \mathbb{R}^{m\times p}\to \mathbb{R}$ is a order-4 tensor, whose components are expressed as: $\frac{\partial ^2f(X)}{\partial X_{ij} \partial X_{kl} }$, where $1\le i,k\le m$ and $1\le j,l\le p$, its Frobenius norm is defined by
\begin{equation}\label{tensorFnorm}
	\|\nabla^2 f(X)\|_F = \sqrt{\underset{i,k=1}{\overset{m}{\sum}}\underset{j,l=1}{\overset{p}{\sum}}|\frac{\partial ^2f(X)}{\partial X_{ij} \partial X_{kl} }|^2} .
\end{equation}
Multiplication of a order-4 tensor and a matrix $V\in \R^{m\times p}$ is defined by
$$[\nabla^2 f(X) V]_{ij} = \underset{k,l} {\sum}\frac{\partial ^2f(X)}{\partial X_{ij} \partial X_{kl} } V_{kl}. $$
We first state two assumptions for the following analysis.
\begin{assumption}\label{Ass2.1}
	The only vector $z\in \R^{|\mathcal{T}_*|}$ satisfying $z\geq 0$ and $A^T_{\mathcal{T}_*} z=0$ is $z=0$.
\end{assumption} 
\begin{assumption}\label{Ass2.2}
	The matrices $A_i$ are linearly independent for $ i \in \mathcal{T}_*$. 
\end{assumption} 
It is clear that Assumption \ref{Ass2.2} implies Assumption \ref{Ass2.1}. By applying the properties of $\|(\cdot)_+\|_0$ and Assumption \ref{Ass2.1}, \ref{Ass2.2}, we have the following propositions and all the proofs are given in appendix. 

\begin{proposition}\label{limsubdiff}
	Suppose {\rm Assumption \ref{Ass2.1}} holds at some $X^*$. Then  $\varphi$ is regular at $X^*$ and
	\begin{align*}
		\partial \varphi(X^*) =\partial^p \varphi(X^*) = \hat{\partial} \varphi(X^*) = \partial^{\infty}\varphi(X^*)=A^T \partial \|F(X^*)_+\|_0
	\end{align*}
	and
	\begin{align*}
		{\rm d}\varphi(X^*)(H) 
		=\left\{\begin{array}{clc}
			0, &H \in {\rm dom}~{\rm d}\varphi(X^*),\\ 
			+\infty, &{\rm otherwise},
		\end{array}
		\right.
	\end{align*}
	where ${\rm dom}~{\rm d}\varphi(X^* ) = \{H\in  \mathbb{R}^{m \times p}~|~\langle A_i,H\rangle \leq 0,~i \in \mathcal{T}_*;~\langle A_i,H\rangle \in \mathbb{R},~i \notin \mathcal{T}_*\}$.
\end{proposition}

\begin{proposition}(critical cone)\label{prop-critica-varphi}
	Suppose {\rm Assumption \ref{Ass2.1}} holds at some $X^*$. Then for any $V\in \partial \varphi(X^*)\in \R^{m\times p}$ and  $z^* \in \Lambda(X^*, V) = \{z \in \partial\|F(X^*)_+\|_0~|~A^Tz = V\}$, the critical cone of $\varphi$ at $X^*$ is  
	\begin{align}\label{critical-varphi}
		K_{\varphi}(X^*,V) 
		=\left\{D  \in \R^{m\times p}~\bigg|~\begin{aligned}
			&z^*_i\langle  A_i,D \rangle=0,~ \langle  A_i,D \rangle\leq 0,\\
			&z^*_i\geq 0,~~i\in \mathcal{T}_*
		\end{aligned}\right\}.
	\end{align}
\end{proposition}

\begin{proposition} \label{Prop 3.6}
	Suppose {\rm Assumption \ref{Ass2.1}} holds at some $X^* \in \mathbb{R}^{m\times p}$, for $H\in {\rm dom}~{\rm d}\varphi(X^*)$, i.e., $AH \in {\rm dom}~{\rm d} g(F(X^*))$, one has that 
	\begin{align*}
		&{\rm d^2}\varphi(X^*)(H,W)\\
		=&{\rm d^2}g(F(X^*))(AH,AW)\\
		=&\left\{
		\begin{array}{clc}
			0, &\begin{aligned}
				&\langle A_i,W \rangle\leq 0, \langle A_i,H \rangle =  0, i \in \mathcal{T}_*;\\
				&\langle A_i,W \rangle\in \mathbb{R},i \notin \mathcal{T}_*,
			\end{aligned}\\ 
			+\infty, &{\rm otherwise},
		\end{array}
		\right.
	\end{align*}
	and	$\varphi$ is parabolically epi-differentiable at $X^*$ for $H$.
\end{proposition}

\begin{proposition}\label{second-derivative}
	Suppose {\rm Assumption \ref{Ass2.1}} holds at some $X^* \in \R^{m\times p}$. Then for any $V\in \partial \varphi(X^*)\in \R^{m\times p} $ and $H \in  K_{\varphi}(X^*, V)$, we have that  
	\begin{align*} 
		{\rm d}^2\varphi(X^*, V)(H)
		=\sup_{z\in \Lambda(X^*,V)}{\rm d}^2g(F(X^*),z)(AH) 
	\end{align*}
	and there exists some $\hat{z} \in \Lambda(X^*,V)$ such that
	\begin{align*} 
		{\rm d}^2\varphi(X^*, V)(H)
		&={\rm d}^2g(F(X^*),\hat{z})(AH) \\
		&= \left\{
		\begin{array}{clc}
			0,&~ H \in K_{\varphi}(X^*,V),\\
			+\infty,&~{\rm otherwise. } 
		\end{array}
		\right.
	\end{align*}
\end{proposition}

At the end of this section, we  define a P-stationary point for problem \eqref{Prim}, which is instrumental in  establishing the P-stationary equation and designing the optimization algorithm.
\begin{definition}
	A point $X^*$ is called a  P-stationary point of \eqref{Prim}  if there exists a constant $ \alpha>0$ and a point $z^*\in \R^n$, such that
\end{definition}
\begin{equation}\label{Proximal}
	\left\{\begin{aligned}
		0 &= \nabla f(X^*)+ A^Tz^*,	\\
		F(X^*) &\in {\rm prox}_{\alpha \lambda \|(\cdot)_+\|_0}(F(X^*) + \alpha z^*).
	\end{aligned}
	\right.
\end{equation}
We also say $(X^*,z^*)$ is a  P-stationary point of problem \eqref{Prim} if \eqref{Proximal} holds. As shown in \cite{ZPXQ},  the proximal operator of $\|(\cdot)_+\|_0$ admits a closed-form expression:
$$[{\rm prox}_{\alpha \lambda \|(\cdot)_+\|_0}(z)]_i=\left\{
\begin{aligned}
	&0,&z_i \in (0,\sqrt{2\lambda \alpha}),\\
	&\{0,z_i\},&z_i \in \{0,\sqrt{2\lambda \alpha}\},\\
	&z_i,&z_i \in(-\infty,0)\cup (\sqrt{2\lambda \alpha},\infty),
\end{aligned}
\right.$$
where $ i=1,2,\cdots,n.$ By the definition of  proximal  operator, the condition $F(X^*) \in {\rm prox}_{\alpha \lambda \|(\cdot)_+\|_0}(F(X^*)+\alpha z^*)$ implies that 
\begin{equation}\label{Proximal2}
	\left\{\begin{aligned}
		&z_i^* \in [0,\sqrt{\frac{2\lambda}{\alpha}}],~ & F_i(X^*) =0,\\
		& F_i(X^*)\in(-\infty,0]\cup [(\sqrt{2\lambda \alpha},\infty)~&  z^*_i = 0.
	\end{aligned}
	\right.
\end{equation}

\subsection{The optimality conditions}\noindent
In this part, we derive the optimality conditions of problem \eqref{Prim} and investigate  the relationship between the P-stationary point and  the local minimizer of problem \eqref{Prim}. As a consequence of \cite[Proposition 3.99]{BS} and Proposition \ref{second-derivative}, we obtain the following necessary and sufficient conditions. 

\begin{theorem}(First-order conditions for optimality)Suppose {\rm Assumption \ref{Ass2.1}} holds at $X^*$. Then the following assertions hold: \\ 
$(\romannumeral 1)$ If $X^*$ is a local minimizer of problem \eqref{Prim}, then 
	\begin{equation*}\label{04}
		\langle \nabla f(X^*), D\rangle \geq 0,
		~\forall~  D \in \left\{D\in  \mathbb{R}^{m\times p}~\bigg|~\begin{aligned}
			&\langle D,A_i\rangle \leq 0,~i \in \mathcal{T}_*;\\
			&\langle D,A_i\rangle \in \mathbb{R},~i \notin \mathcal{T}_*
		\end{aligned}\right\} .
	\end{equation*}
$(\romannumeral 2)$ If  the condition 
	\begin{equation*}\label{04}
		\langle \nabla f(X^*), D\rangle > 0,~\forall~ D \in \left\{D\in  \mathbb{R}^{m\times p}\bigg|\begin{aligned}
			&\langle D,A_i\rangle \leq 0,~i \in \mathcal{T}_*;\\
			&\langle D,A_i\rangle \in \mathbb{R},~i \notin \mathcal{T}_*
		\end{aligned}\right\} \setminus\{0\},
	\end{equation*}
	then $X^*$ is a local minimizer of problem \eqref{Prim}.
\end{theorem}

\begin{theorem}(Second-order  conditions for optimality)\label{Second-condition}
	Suppose {\rm Assumption \ref{Ass2.1}} holds at $X^*$. Then the following assertions hold:\\
	$(\romannumeral 1)$If $X^*$ is a local minimizer of problem \eqref{Prim}, then $ 0 \in \partial(f+\lambda\varphi)(X^*)$,~i.e., $ -\nabla f(X^*)\in \lambda\partial \varphi(X^*) =  \partial\varphi(X^*)$, there exists $ z^*\in \varLambda(X^*,-\nabla f(X^*))$ such that 
	\begin{equation*}\label{04}
		\langle D, \nabla^2 f(X^*)D  \rangle \geq 0,
		~\forall~  D \in\left\{\begin{aligned} 
			&D  \in \R^{m\times p}~|~z^*_i \langle D,A_i\rangle =0,\\
			&\langle D,A_i \rangle \leq 0,~z^*_i\geq 0,~~i\in \mathcal{T}_*
		\end{aligned}\right\}.
	\end{equation*}	
$(\romannumeral 2)$ If $ 0 \in \partial(f+\lambda\varphi)(X^*)$  and there exists $ z^*\in \varLambda(X^*,-\nabla f(X^*))$ such that  
	$$\langle D, \nabla^2 f(X^*)D  \rangle > 0,~\forall~D \in \left\{\begin{aligned} 
		&D  \in \R^{m\times p}~|~z^*_i \langle D,A_i\rangle =0,\\
		&\langle D,A_i \rangle \leq 0,~z^*_i\geq 0,~i\in \mathcal{T}_*
	\end{aligned}\right\}\setminus \{0\}.$$
	Then, $X^*$ is a local minimizer of problem \eqref{Prim} and the second order growth condition holds at $X^*$, i.e., there exist $\kappa>0$ and $\delta >0$ such that for any $X \in B(X^*,\delta)$, 
	\begin{equation*}
		f(X)+\lambda\|F(X)_+\|_0 \ge f(X^*)+\lambda\|F(X^*)_+\|_0 + \kappa \|X-X^*\|_F^2.
	\end{equation*} 
\end{theorem}

To apply the above optimality conditions to problem \eqref{QSSVM0/1}, we first introduce some notions 
\begin{equation}\label{f(w,b,c)}
	f(W, b,c) =\sum_{i=1}^n \frac{1}{2}\|Wx^i + b\|^2,
\end{equation}
$$F(W,b,c) =  \begin{pmatrix}
	1-y_1(\frac{1}{2} {x^1}^T Wx^1 + b^T x^1+c)\\
	1-y_2(\frac{1}{2} {x^2}^T Wx^2 + b^T x^2+c)\\
	\vdots\\
	1-y_n(\frac{1}{2} {x^n}^T Wx^n + b^T x^n+c)
\end{pmatrix},$$
and 
$$\varphi(W,b,c)  = \underset{i=1}{\overset{n}{\sum}}\ell_{0/1}(1-y_i(\frac{1}{2} {x^i}^T Wx^i + b^T x^i+c)).$$

The gradient of $f$ and $F$ at $(W^*,b^*,c^*)$ is, respectively
$$\nabla f(W^*,b^*,c^*) =\overset{n}{\underset{i=1}{\sum}}  \begin{pmatrix}
	W^*x^i{x^i}^T+b^*{x^i}^T\\
	W^*x^i +b^*\\
	0
\end{pmatrix},$$
$$\nabla F(W^*,b^*,c^*) = -\begin{pmatrix}
	\frac{1}{2} y_1x^1{x^1}^T& y_1{x^1}^T &y_1\\
	\frac{1}{2} y_2x^2{x^2}^T& y_2{x^2}^T &y_2\\
	\vdots &\vdots & \vdots\\
	\frac{1}{2} y_nx^n{x^n}^T& y_n{x^n}^T &y_n\\
\end{pmatrix}.$$

\begin{assumption}\label{Ass2.3}
The only vector $z\in \R^{|\tilde{\mathcal{T}}_*|}$ satisfying $z\geq 0$ and $\nabla F(W^*,b^*,c^*) ^T_{\tilde{\mathcal{T}}_*} z=\begin{pmatrix}
0_{m\times m}\\ 0_m\\0
\end{pmatrix}$ is $z=0$, where $\tilde{\mathcal{T}}_* := \{	i~|~y_i(\frac{1}{2} {x^i}^T W^*x^i + {b^*}^T x^i+c^*)=1\}$.
\end{assumption} 
Then we have the following theorems.
\begin{theorem}(First-order conditions for optimality)
	Suppose {\rm Assumption \ref{Ass2.3}} holds at $(W^*,b^*,c^*)^*$. Then the following assertions hold:\\		
	$(\romannumeral 1)$ If $(W^*,b^*,c^*)$ is a local minimizer of problem \eqref{QSSVM0/1}, then 
	\begin{equation*}\label{04}
		\overset{n}{\underset{i=1}{\sum}} \langle W^*x^i{x^i}^T +b^*{x^i}^T, D_1\rangle +\langle W^*x^i +b^*, d_2 \rangle \geq 0,
	\end{equation*}
	for any $ (D_1, d_2, d_3) \in \{(D_1, d_2, d_3)\in  \mathbb{R}^{m\times m}\times \mathbb{R}^{m}\times \R~|~y_i(\frac{1}{2} {x^i}^T D_1x^i + d_2^Tx^i+d_3) \ge 0,~i \in \mathcal{T}_*;~y_i(\frac{1}{2} {x^i}^T D_1x^i + d_2^Tx^i+d_3) \in \mathbb{R},~i \notin \mathcal{T}_*\} .$\\
	$(\romannumeral 2)$If  the condition 
	\begin{equation*}\label{04}
		\overset{n}{\underset{i=1}{\sum}} \langle W^*x^i{x^i}^T +b^*{x^i}^T, D_1\rangle +\langle W^*x^i +b^*, d_2 \rangle >0,
	\end{equation*}
	for any $ (D_1, d_2, d_3) \in \{(D_1, d_2, d_3)\in  \mathbb{R}^{m\times m}\times \mathbb{R}^{m}\times \R~|~y_i(\frac{1}{2} {x^i}^T D_1x^i + d_2^Tx^i+d_3) \ge 0,~i \in \mathcal{T}_*;~y_i(\frac{1}{2} {x^i}^T D_1x^i + d_2^Tx^i+d_3) \in \mathbb{R},~i \notin \mathcal{T}_*\}\setminus \{0\}$, then $(W^*,b^*,c^*)$ is a local minimizer of problem \eqref{QSSVM0/1}.
\end{theorem}

\begin{theorem}(Second-order conditions for optimality)
		Suppose {\rm Assumption \ref{Ass2.3}} holds at $(W^*,b^*,c^*)^*$. Then the following assertions hold:\\
	$(\romannumeral 1)$ If $(W^*,b^*,c^*)$ is a local minimizer of problem \eqref{QSSVM0/1}, then $0\in \partial(f+\lambda\varphi)(W^*,b^*,c^*)$, i.e., $ -\nabla f(W^*,b^*,c^*)\in \lambda\partial \varphi(W^*,b^*,c^*)=\partial\varphi(W^*,b^*,c^*)$, there exists $z^*\in \varLambda((W^*,b^*,c^*),-\nabla f(W^*,b^*,c^*)) = \{z \in \partial\|F(W^*,b^*,c^*)_+\|_0~|~\nabla F(W^*,b^*,c^*)^Tz =-\nabla f(W^*,b^*,c^*)\}$ such that 
	\begin{equation*}\label{04}
		\overset{n}{\underset{i=1}{\sum}} \|D_1x^i+d_2\|^2  \ge 0,
	\end{equation*}
	for any $ (D_1, d_2, d_3) \in \{(D_1, d_2, d_3)\in  \mathbb{R}^{m\times m}\times \mathbb{R}^{m}\times \R~|~\frac{1}{2} {x^i}^T D_1x^i + d_2^Tx^i+d_3 = 0,~z^*_i \ge 0; ~\frac{1}{2} {x^i}^T D_1x^i + d_2^Tx^i+d_3 \neq 0,~z^*_i = 0,~i \in \mathcal{T}_*\}$\\
	$(\romannumeral 2)$ If $ 0 \in \partial(f+\lambda\varphi)(W^*,b^*,c^*)$  and there exists $ z^*\in \varLambda((W^*,b^*,c^*),-\nabla f(W^*,b^*,c^*))$ such that  
	$$	\overset{n}{\underset{i=1}{\sum}} \|D_1x^i+d_2\|^2 > 0,$$
	for any $ (D_1, d_2, d_3) \in \{(D_1, d_2, d_3)\in  \mathbb{R}^{m\times m}\times \mathbb{R}^{m}\times \R~|~\frac{1}{2} {x^i}^T D_1x^i + d_2^Tx^i+d_3 = 0,~z^*_i \ge 0; ~\frac{1}{2} {x^i}^T D_1x^i + d_2^Tx^i+d_3 \neq 0,~z^*_i = 0,~i \in \mathcal{T}_*\}\setminus \{0\}$.
	Then, $(W^*,b^*,c^*)$ is a local minimizer of problem \eqref{QSSVM0/1} and the second order growth condition holds at $(W^*,b^*,c^*)$, i.e., there exist $\kappa>0$ and $\delta >0$ such that for any $(W^*,b^*,c^*) \in B((W^*,b^*,c^*),\delta)$, 
	\begin{equation*}
		\begin{aligned}
			&f(W,b,c)+\lambda\|F(W,b,c)_+\|_0 \\
			\ge& f(W^*,b^*,c^*)+\lambda\|F(W^*,b^*,c^*)_+\|_0\\
			&+ \kappa (\|W-W^*\|_F^2+\|(b,c) -(b^*,c^*)\|^2).
		\end{aligned}
	\end{equation*} 
\end{theorem}

\begin{theorem}\label{stationary-equation}
	The following assertions hold for problem \eqref{Prim}\\
	$(\romannumeral 1)$  If $X^*$ is a local minimizer, then there exists  $z^*$ such that $(X^*, z^*)$ is a P-stationary point for any $0<\alpha <\alpha_*:=\min\{\alpha_1, \alpha_2\}$, where 
	\begin{equation}\label{g_i(X^*)}
		\begin{aligned}
			&\alpha_1 = \left\{\begin{aligned}
				& + \infty, ~ &F_i(X^*)\leq 0,\\
				&\min\{\frac{F_i(X^*)^2}{2\lambda}:~F_i(X^*)>0\},~&{\rm otherwise,}
			\end{aligned}
			\right.\\
			&\alpha_2 = \left\{\begin{aligned}
				& + \infty, ~ &z_i^*\leq 0,\\
				&\min\{\frac{2\lambda}{(z_i^*)^2}:~z_i^*>0\},~&{\rm otherwise}.
			\end{aligned}
			\right.
		\end{aligned}
	\end{equation}
	$(\romannumeral 2)$A point $X^*$ is a local minimizer of problem \eqref{Prim} if  $(X^*,z^*)$  is a P-stationary point with $\alpha >0$.			 
\end{theorem}

\section{ Newton method }
In this section, we shall establish a  Newton method for the optimization problem \eqref{Prim} and analyze its local quadratic convergence under some reasonable assumptions and the second-order sufficient condition. Our approach is going along the idea  of \cite{ZPXQ}.  We first present P-stationary equation for \eqref{Prim}, then based on it construct a Newton method, and prove local quadratic convergence of this method. Our result extend the corresponding results of \cite{ZPXQ} to matrix case. 

\subsection{Stationary equation}
In this part, we derive a stationary equation, which plays a significant role in algorithmic design. Firstly,  we define the following index sets for a point $(X, z)$:
\begin{align}\label{index}
	\mathcal{T}_o &: = \{F_i(X) = 0,~ \alpha z_i \in\{0, \sqrt{2\alpha \lambda}\}\},\notag\\
	\mathcal{T}_1 &: = \{F_i(X) + \alpha z_i \in(0, \sqrt{2\alpha \lambda})\},\notag\\
	\mathcal{T}_2 &: = \{F_i(X) + \alpha z_i \in\{0, \sqrt{2\alpha \lambda}\}\},\\
	\mathcal{T}_3 &: = \{F_i(X) + \alpha z_i \in(-\infty, 0)\cup ( \sqrt{2\alpha \lambda}, +\infty)\},\notag
\end{align}
similarly, we defined the index sets $\mathcal{T}_o^*$,  $\mathcal{T}_1^*$, $\mathcal{T}_2^*$, $\mathcal{T}_3^*$ based on $(X^*,z^*)$. For a fixed index set $\mathcal{T}$, we define the  stationary equation 
\begin{equation}\label{stationary equation }
	\Psi(X,z;\mathcal{T}) = 
	\begin{pmatrix}
		\nabla f(X)+ A^T_\mathcal{T}z_\mathcal{T}\\
		F(X)_\mathcal{T}\\
		z_{\bar{\mathcal{T}}}
	\end{pmatrix} = 0
\end{equation}
and $\|\Psi(X,z;\mathcal{T})\|_F = (\|\nabla f(X)+ A^T_\mathcal{T}z_\mathcal{T}\|_F^2 + \|(F(X)_\mathcal{T},z_{\bar{\mathcal{T}}})\|^2)^{\frac{1}{2}}$. It follows from \cite[Theorem 4.3]{ZPXQ} that $(X^*,z^*)$ is a P-stationary point of problem \eqref{Prim} with $\alpha \ge 0$  if and only if  $\Psi(X^*,z^*;\mathcal{T}_*)=0$ and $\mathcal{T}_*=  \mathcal{T}_o^* \cup\mathcal{T}_1^*$.

\subsection{Algorithmic design}
Based on  the conclusion in Theorem \ref{stationary-equation}, we are going to establish the Newton method  to solve the stationary equation \eqref{stationary equation }. A noteworthy issue is that we can not obtain the index $\mathcal{T}_*$ if the point $X^*$ is unknown. Therefore, we update index set and the iteration points alternatively. Let $(X^k,z^k)$ be the $k$-th iteration point and its associated  index set is defined as
\begin{equation}\label{J_k}
	\mathcal{T}_k =\mathcal{T}_o^k\cup \mathcal{T}_1^k .
\end{equation}
Let $D^k=(D^k_{X}, D^k_{\mathcal{T}_k}, D^k_{\bar{\mathcal{T}}_k})\in \R^{m\times p} \times \R^{|\mathcal{T}_k|}\times \R^{|\bar{\mathcal{T}}_k|}$ be the Newton direction, which  satisfies the following equation:
\begin{align}\label{Newton equation}
	\nabla \Psi(X^k, z^k; \mathcal{T}_k) D^k = -\Psi(X^k, z^k; \mathcal{T}_k).
\end{align}
Here, $\nabla \Psi(X^k, z^k; \mathcal{T}_k)$ is the Jacobian matrix of $\Psi(X^k, z^k; \mathcal{T}_k)$  with respect to $(X^k, z^k)$ and it is a order-4 tensor.  Instead of solving the above equation \eqref{Newton equation}, we try to solve its perturbed
version, we replace $\nabla \Psi(X^k, z^k; \mathcal{T}_k) $  with $\nabla \Psi_{{\gamma}_k}(X^k, z^k; \mathcal{T}_k)$, then the Newton direction $(D^k_{X}, D^k_{\mathcal{T}_k}, D^k_{\bar{\mathcal{T}}_k})$ satisfies
\begin{equation}\label{Newton-direction}
	\left\{\begin{aligned}
		\nabla^2 f(X^k)D^k_{X} + A^T_{\mathcal{T}_k}D^k_{\mathcal{T}_k} &= -(\nabla f(X)+ A^T_{\mathcal{T}_k}z_{\mathcal{T}_k })\\
		A_{\mathcal{T}_k}D^k_{X} -{\gamma}_k D^k_{\mathcal{T}_k} &= -F(X)_{\mathcal{T}_k}\\
		D^k_{\bar{\mathcal{T}}_k} &= -z_{\bar{\mathcal{T}}_k}
	\end{aligned}
	\right.
\end{equation}		
The rule for updating ${\gamma}_k$ is as follows 
\begin{equation}\label{gamma}
	{\gamma}_k: = \min\{\tau\gamma_{k-1}, \rho\|\Psi(X^k, z^k; \mathcal{T}_k)\|_F\},
\end{equation}
where $\tau\in(0,1)$ and $\rho>0$.	The  Newton method is summarized in Algorithm \ref{alg:example}.
\begin{algorithm}[H]
	\caption{ Newton Method for problem \eqref{Prim}.}\label{alg:example}
	\begin{algorithmic}
		\STATE 
		\STATE {\textbf{Initialization}} $(X^0,z^0)$  and $\gamma_{-1}>0$. Set the parameter $\alpha,~\rho,~\lambda>0$, $\tau \in (0,1)$. Compute $\mathcal{T}_0$ by \eqref{J_k} .
		\STATE \hspace{0.5cm}$ \textbf{for} k = 0,1,\cdots$ \textbf{do} 
		\STATE \hspace{0.5cm}\textbf{If} $\|\Psi(X^k,z^k;\mathcal{T}_k)\|_F>0$, \textbf{then }
		\STATE \hspace{0.5cm}Compute ${\gamma}_k$  by \eqref{gamma} and Newton direction $D^k$  by \eqref{Newton-direction}.
		\STATE \hspace{0.5cm}Set $(X^{k+1},z^{k+1})=(X^k, z^k)+D^k$ . 
		\STATE \hspace{0.5cm}Update $\mathcal{T}_{k+1}$  by \eqref{J_k} and ${\gamma}_{k+1}$  by \eqref{gamma}, set $k := k+1$.
		\STATE \hspace{0.5cm}\textbf{end if}
		\STATE \hspace{0.5cm}\textbf{return} $(X^k,z^k)$.
	\end{algorithmic}
	\label{alg1}
\end{algorithm}
	\subsection{Local quadratic convergence}
To prove the local quadratic convergence of the proposed method, we make the following assumption.
\begin{assumption}\label{Ass4.1}
	The Hessian matrix $\nabla^2f $ is local Lipschitz continuous around $X^*$ with a constant $L_*>0$. That is,  there exists  $\delta_0^*>0$ such that
	\begin{align*}
		\|\nabla^2f (X^1)- \nabla^2f (X^2)\|_F\leq L_*\|X^1-X^2\|_F, &\\
		  \forall~ X^1, X^2&\in B(X^*,\delta_0^*).
	\end{align*}
\end{assumption}
\begin{remark}
	It is clear that the function $f$ denoted as \eqref{f(w,b,c)} in  problem \eqref{QSSVM0/1} satisfies Assumption \ref{Ass4.1} around $(W^*,b^*,c^*)$ with $L_* =0$.
\end{remark} 

The proofs of the following lemmas and theorem follow the similar procedure in \cite{ZPXQ}, and all the proofs are given in appendix.  Firstly, we will analyze the properties of the P-stationary point equation in a neighborhood of $(X^*, z^*)$. For convenience, denoted by $ u = F(X)$ and $ u^* = F(X^*)$.
\begin{lemma}\label {4.1}
	Let $(X^*, z^*)$ be a P-stationary point with $0<\alpha<\alpha_*:=\min\{\alpha_1,\alpha_2\}$ of problem \eqref{Prim} and let $\mathcal{T}_* = \mathcal{T}_o^*\cup \mathcal{T}_1^*$. Then there exists  $\delta_1^*>0$ such that for any $(X,z) \in B((X^*,z^*),\delta^*_1)$ with its associated indices $\mathcal{T}_1$ and $\mathcal{T}_o$, it holds that
	$$\Psi(X^*,z^*;\mathcal{T})=0~ {\rm and}~ \mathcal{T}=\mathcal{T}_o\cup \mathcal{T}_1\subseteq \mathcal{T}_*.$$
\end{lemma}

In order to calculate the Newton direction, we introduce the vectorization  of $\Psi(X^k, z^k; \mathcal{T}_k) $ :
$$\tilde{\Psi}(X^k, z^k; \mathcal{T}_k) = \begin{pmatrix}
	{\rm vec}(\nabla f(X^k))+ {\rm vec}(\underset{i\in {\mathcal{T}_k} }{\sum}z_i^k A_i)\\
	F(X^k)_{\mathcal{T}_k}\\
	z_{\bar{\mathcal{T}_k}}^k
\end{pmatrix}\in \R^{mp+n}.$$
Moreover, It follows from \cite{golub2013matrix} that the tensor $\nabla \Psi_{\gamma_k}(X^k,z^k;\mathcal{T}_k)$ can be written as the following matrix form by vectorizing  matrix:
\begin{align*}
	\nabla \tilde{\Psi}_{\gamma_k}(X^k,z^k;\mathcal{T}_k) &:= \begin{pmatrix}
		\frac{\partial^2 f(X^k)}{\partial{\rm vec}(X^k)\partial{\rm vec}(X^k)^T}&  B_{\mathcal{T}_k}^T& 0\\
		B_{\mathcal{T}_k}& -{\gamma_k} I_{|\mathcal{T}_k|} &0\\
		0&0&I_{|\bar{\mathcal{T}}_k|}
	\end{pmatrix}\\
	&\in \R^{(mp+n) \times (mp+n)},
\end{align*}
where $B_{\mathcal{T}_k} = [{\rm vec}(A_i), \cdots, {\rm vec}(A_j)]_{i,j \in{\mathcal{T}_k}}$ $ \in \R^{mp \times |\mathcal{T}_k|}$. From \eqref{tensorFnorm}, we have the Frobenius norm of  the two forms  are equal, i.e.,
$$\|\nabla \Psi_{\gamma_k}(X^k,z^k;\mathcal{T}_k)\|_F = \|\nabla \tilde{\Psi}_{\gamma_k}(X^k,z^k;\mathcal{T}_k)\|_F.$$ 
It is clear that equation \eqref{Newton-direction} is equivalent to
$$ \nabla \tilde{\Psi}_{\gamma_k}(X^k,z^k;\mathcal{T}_k){\rm vec} (D^k)  = -\tilde{\Psi}(X^k, z^k; \mathcal{T}_k).$$		

According to the above equivalence,  we discuss the properties of   $\nabla \tilde{\Psi}_{\gamma_k}(X^k,z^k;\mathcal{T}_k)$ instead of the tensor $\nabla \Psi_{\gamma_k}(X^k,z^k;\mathcal{T}_k)$ in the following part. The next lemma shows the nonsingularity of the Jacobian matrix $\nabla \tilde{\Psi}_{\gamma}(X,z;\mathcal{T})$ over a neighborhood of a P-stationary point $ (X^*,z^*)$. 
\begin{lemma}\label{4.2}
	Let $(X^*,z^*)$ be a P-stationary point with $0<\alpha<\alpha_*$ of problem \eqref{Prim}. Assume {\rm Assumption \ref{Ass2.2}}  and {\rm \ref{Ass4.1}}  hold and the second-order sufficient condition of problem \eqref{Prim} holds. Then there exist two positive constants $c_*$ and $C_*$ such that
	\begin{equation}\label{15}
		C_*\ge \sigma_{\rm max}(\nabla \tilde{\Psi}_{\gamma}(X,z;\mathcal{T}))\ge \sigma_{\rm min}(\nabla \tilde{\Psi}_{\gamma}(X,z;\mathcal{T}))\ge c_*>0,
	\end{equation}
	for any $(X,z) \in B((X^*,z^*), \delta_2^*)$ and any $ 0<\gamma <\frac{c_*}{2}$, where 
	\begin{equation}\label{parameter1}
		C_* := 2{\rm max}\{1, \sigma_{\rm max}(H(\mathcal{T}_*))\}, ~ c_*:=\frac{1}{2} {\rm min}\{1, \underset{\mathcal{T} \subseteq \mathcal{T}_*}{\rm min}\sigma_{\rm min}(H(\mathcal{T}))\}.
	\end{equation}
	\begin{equation}\label{parameter2}
		\delta_2^* := {\rm min}\{\delta^*_1, \frac{c_*}{2L_*}\}, 
		~ H(\mathcal{T}):= \begin{pmatrix}
			\frac{\partial^2 f(X^k)}{\partial{\rm vec}(X^k)\partial{\rm vec}(X^k)^T}&  B^T_{\mathcal{T}}\\
			B_{\mathcal{T}}&  0
		\end{pmatrix}.
	\end{equation}
\end{lemma}

Next, we present the following local quadratic convergence theorem for the  Newton method.
\begin{theorem}\label{convergence}
	Let $(X^*,z^*)$ be a P-stationary point with $0<\alpha<\alpha_*$ for problem \eqref{Prim}. $\delta_2^*$, $c_*$, $C_*$ be given by \eqref{parameter1} and \eqref{parameter2}. Assume {\rm Assumption \ref{Ass2.2}}  and {\rm \ref{Ass4.1}}  hold at $X^*$ and the second-order sufficient condition of problem \eqref{Prim} holds. Let $\{(X^k,z^k)\} $ be the sequence generated by {\rm Algorithm 1} with $0\le \gamma_{-1}\le\frac{c_*}{2}$. If the initial point satisfies $(X^0,z^0) \in  B((X^*,z^*),\delta^*)$,  where
	\begin{equation}\label{delta2}
		\delta^*:={\rm min}\{\delta_2^*,~\frac{c_*}{2(L_*+2\rho C_*)}\},
	\end{equation}
	then the following results hold.\\
	$(\romannumeral 1)$The sequence $\{D^k\}_{k=1}^{\infty}$ is well defined and $\underset{k \to \infty}{\lim } (D^k_{X}, D^k_{\mathcal{T}_k}, D^k_{{\bar{\mathcal{T}}_k}}) =0.$\\
	$(\romannumeral 2)$The whole sequence $\{(X^k,z^k)\}$ convergence to $(X^*,z^*)$ quadratically, namely,
	\begin{align*}
		\|X^{k+1} - X^*\|_F +\|z^{k+1} -z^*\| \le &\frac{L_*+2\rho C_*}{c_*}(\|X^k- X^*\|_F^2 \\
		&+\|z^{k+1} -z^*\|^2).
	\end{align*}
	$(\romannumeral 3)$ The halting condition satisfies 
	$$\|\Psi(X^{k+1},z^{k+1};\mathcal{T}_{k+1})\|_F\le \frac{C_*(L_*+2\rho C_*)}{c_*^3}\|\Psi(X^k,z^k;\mathcal{T}_k)\|_F^2 $$	
	and  {\rm Algorithm 1} achieves $\|\Psi(X^k,z^k;\mathcal{T}_k)\|<\epsilon$ for a given tolerance $\epsilon>0$ when
	\begin{align}\label{epsilon}
		k\ge& {\rm log}_2(2\sqrt{\frac{(L_*+2\rho C_*)C_*^3}{c_*^3}}(\|X^0 - X^*\|_F +\|z^0 -z^*\|))\notag\\
		&-{\rm log}_2\sqrt{\epsilon}.
	\end{align}
\end{theorem}

\begin{remark}
	{\rm Theorem \ref{convergence}} shows that the  Newton method is locally quadratic convergent. It is worth noting that the initial point should be close to the optimal solution and  different initial points may result in different output solutions in our method. In practical applications, we can use some first-order gradient descent methods or convex relaxation methods as an initialization step to obtain an appropriate starting point.
\end{remark}

	\section{Numerical Experiments}
In this section, we conduct extensive numerical experiments to evaluate the performance of Algorithm 1 for problem \eqref{QSSVM0/1}. The experiments are performed using MATLAB (R2024b) on a laptop with 16GB memory and Intel(R) Core(TM) i5-12500H 2.50 GHz CPU. We compare the Newton method for $L_{0/1}$-SQSSVM and several leading solvers for SQSSVM problems, including the ADMM method for QSSVM$_{0/1}$ \cite{Yang}, SVM with polynomial kernel and QSSVM. The evaluation utilizes two types of datasets: artificially generated data and public benchmark data sets.

\subsection{Experiments on artificial data sets}
\begin{figure}[H]
	\centering
	\subfloat[]{\includegraphics[width= 1.7in, height= 1in ]{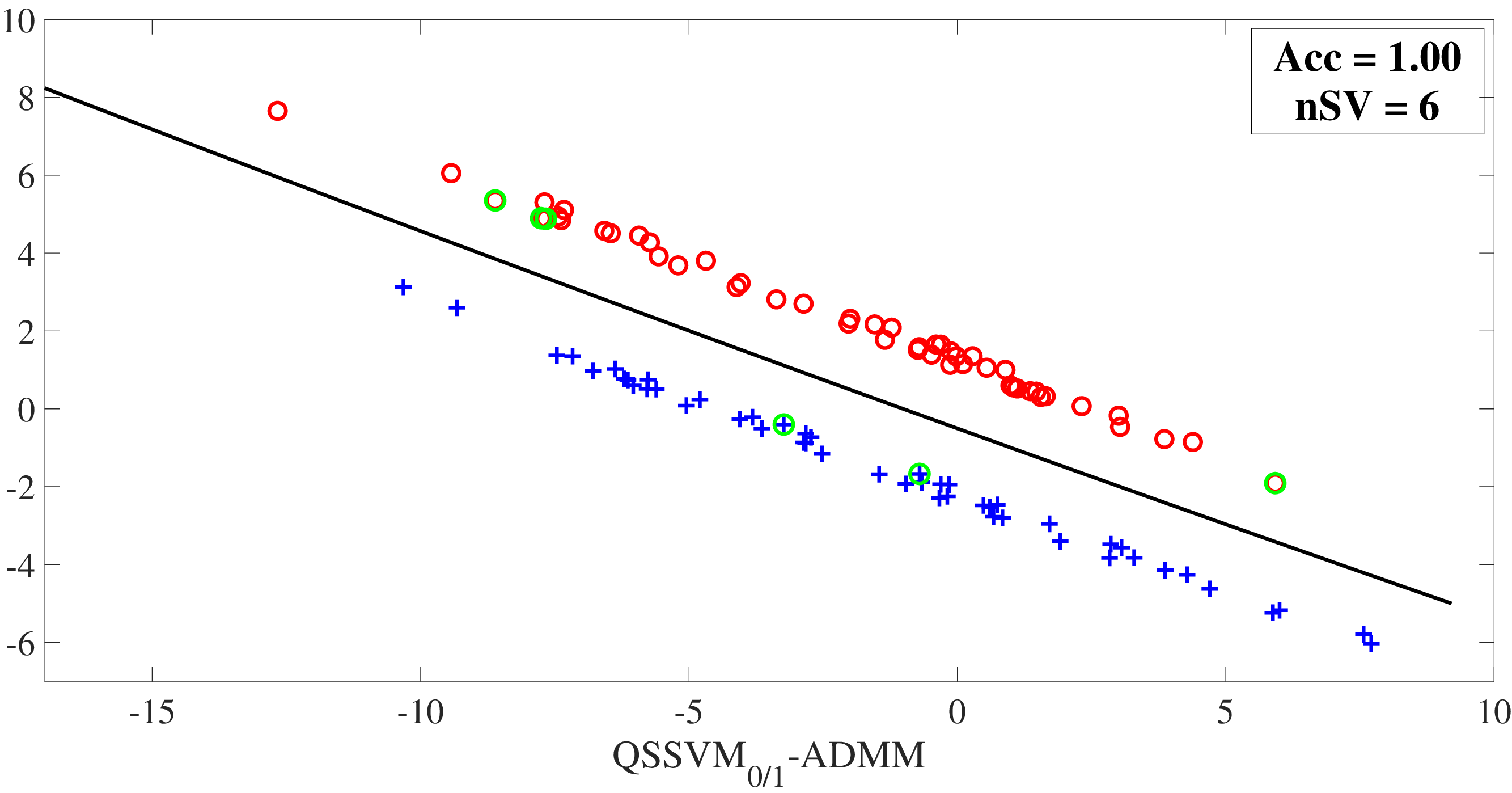}}%}
	\hfil
	\subfloat[]{\includegraphics[width=1.7in, height= 1in]{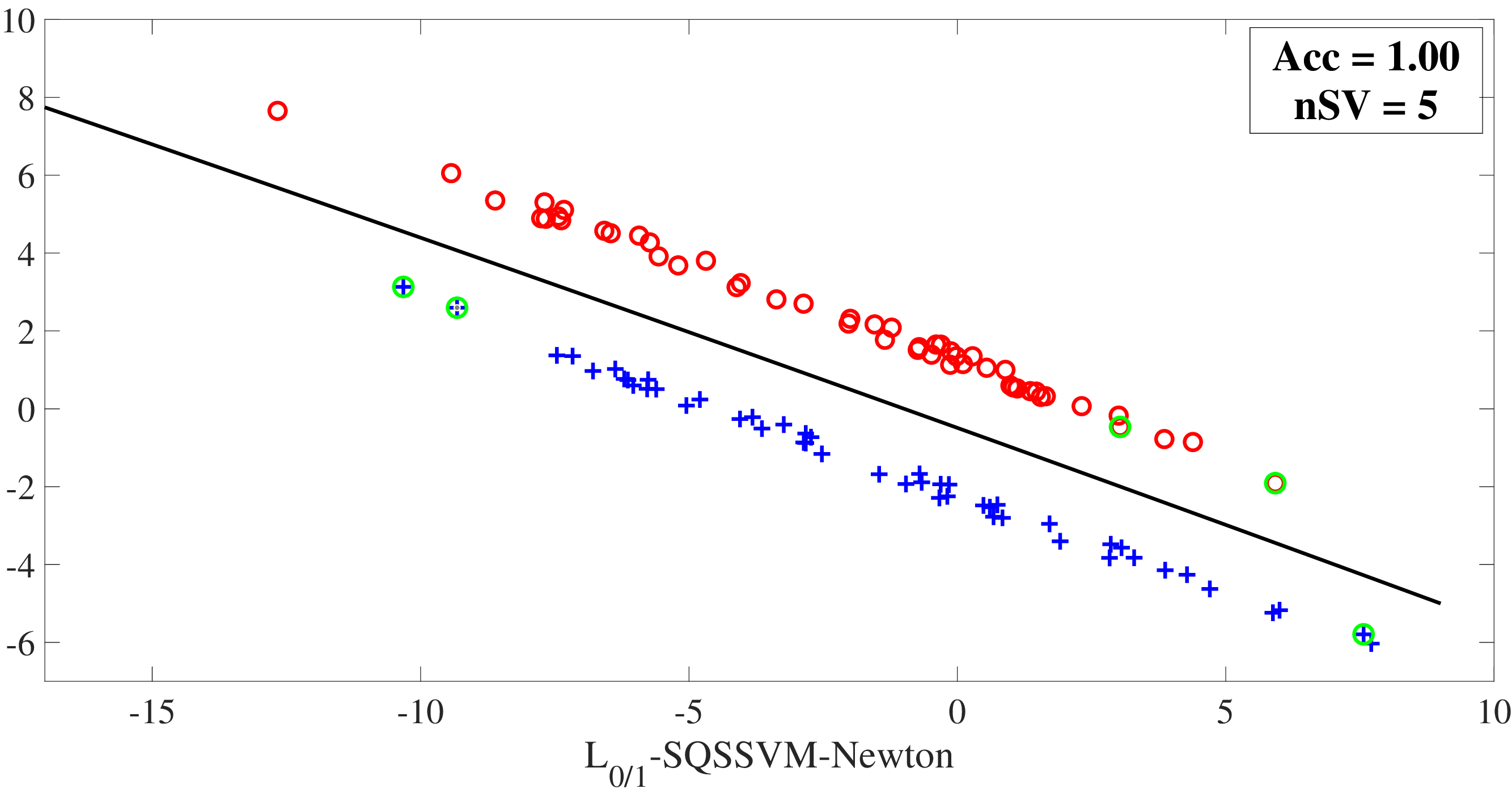}}\\
	\subfloat[]{\includegraphics[width=1.7in, height= 1in]{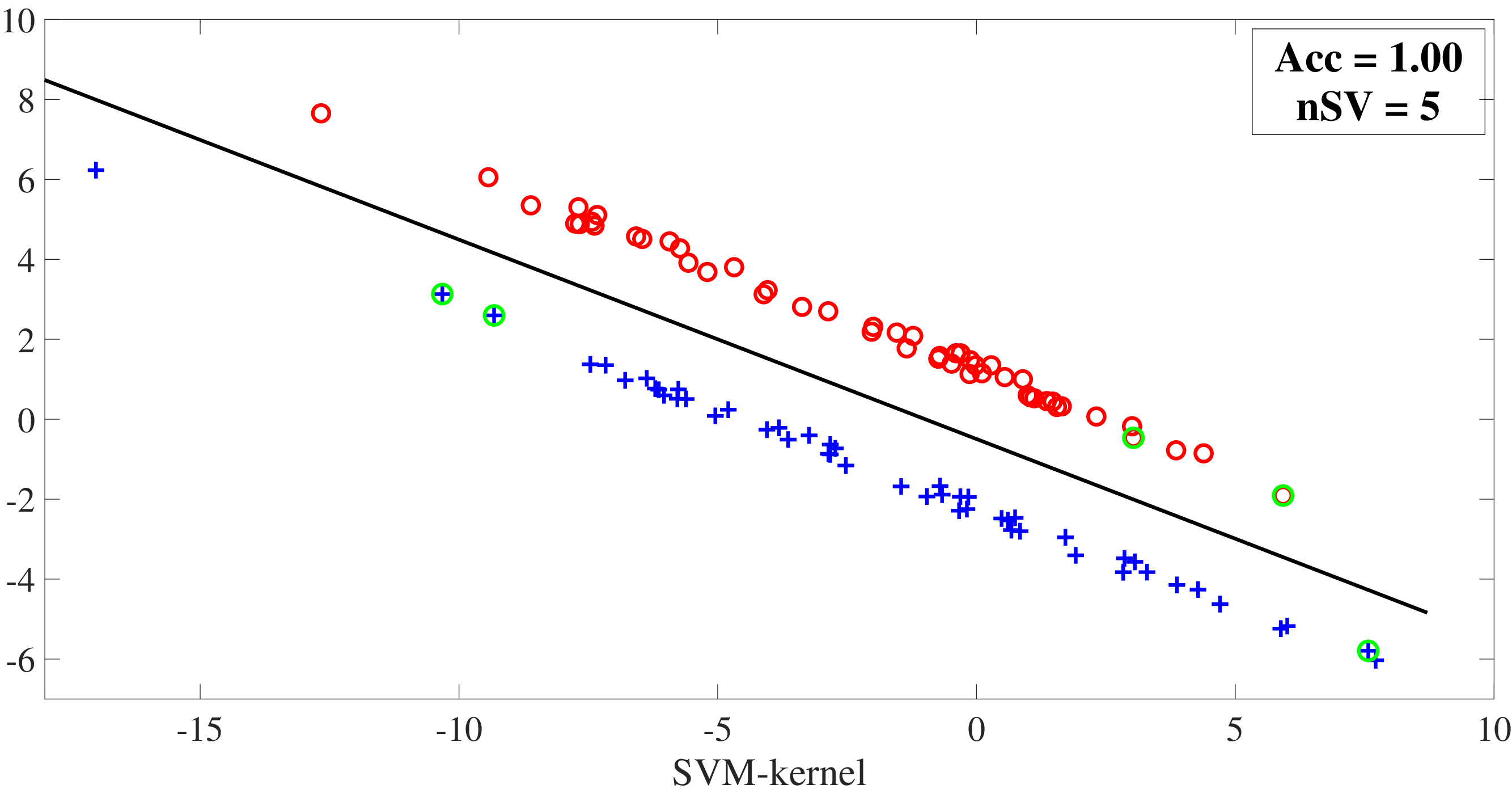}}%
	\hfil
	\subfloat[]{\includegraphics[width=1.7in, height= 1in]{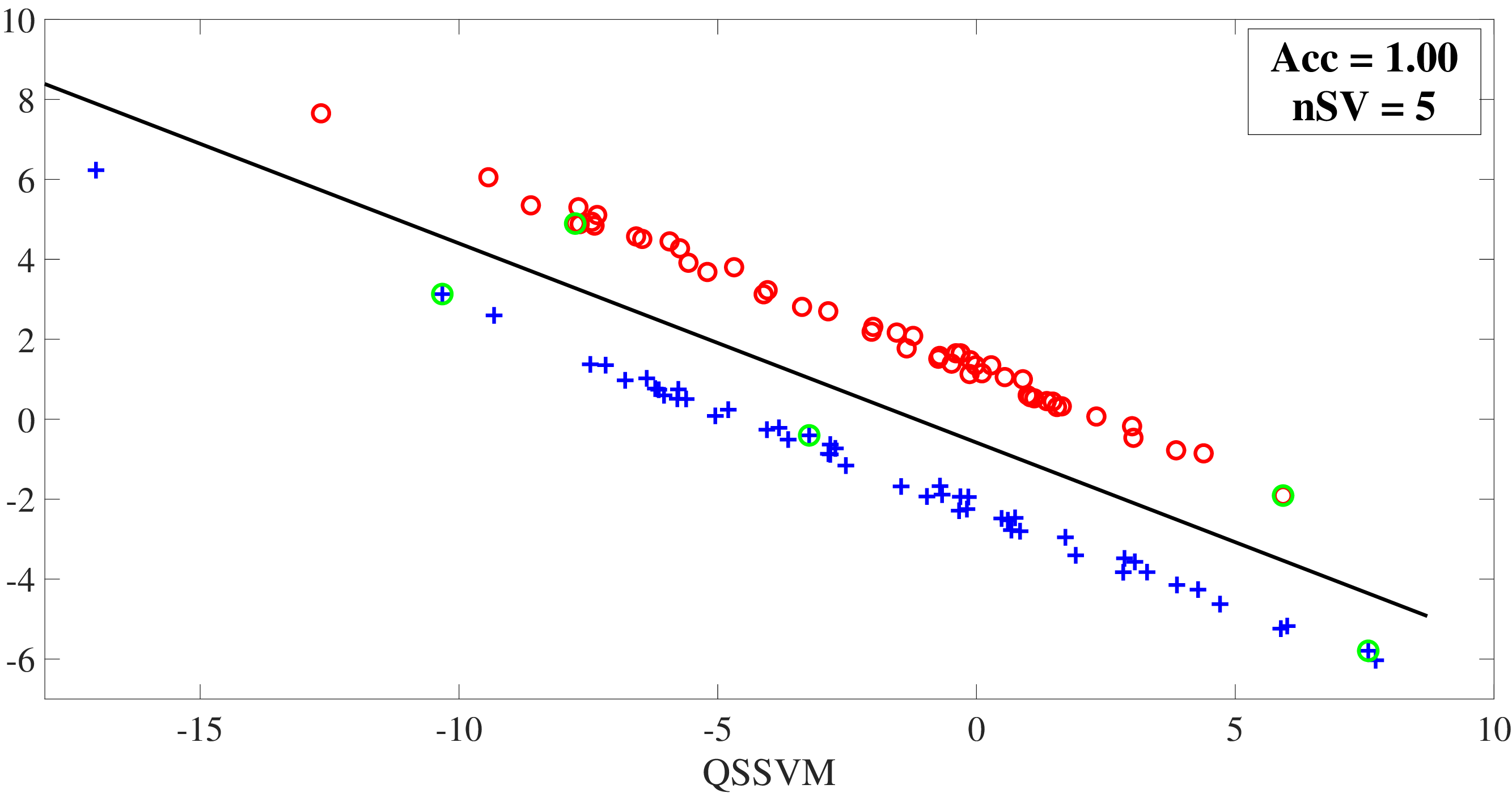}}%
	\caption{Linearly separable data set with 100 points.}
	\label{lindata}
\end{figure}

\begin{figure}[H]
\centering
\subfloat[]{\includegraphics[width=1.7in, height= 1in]{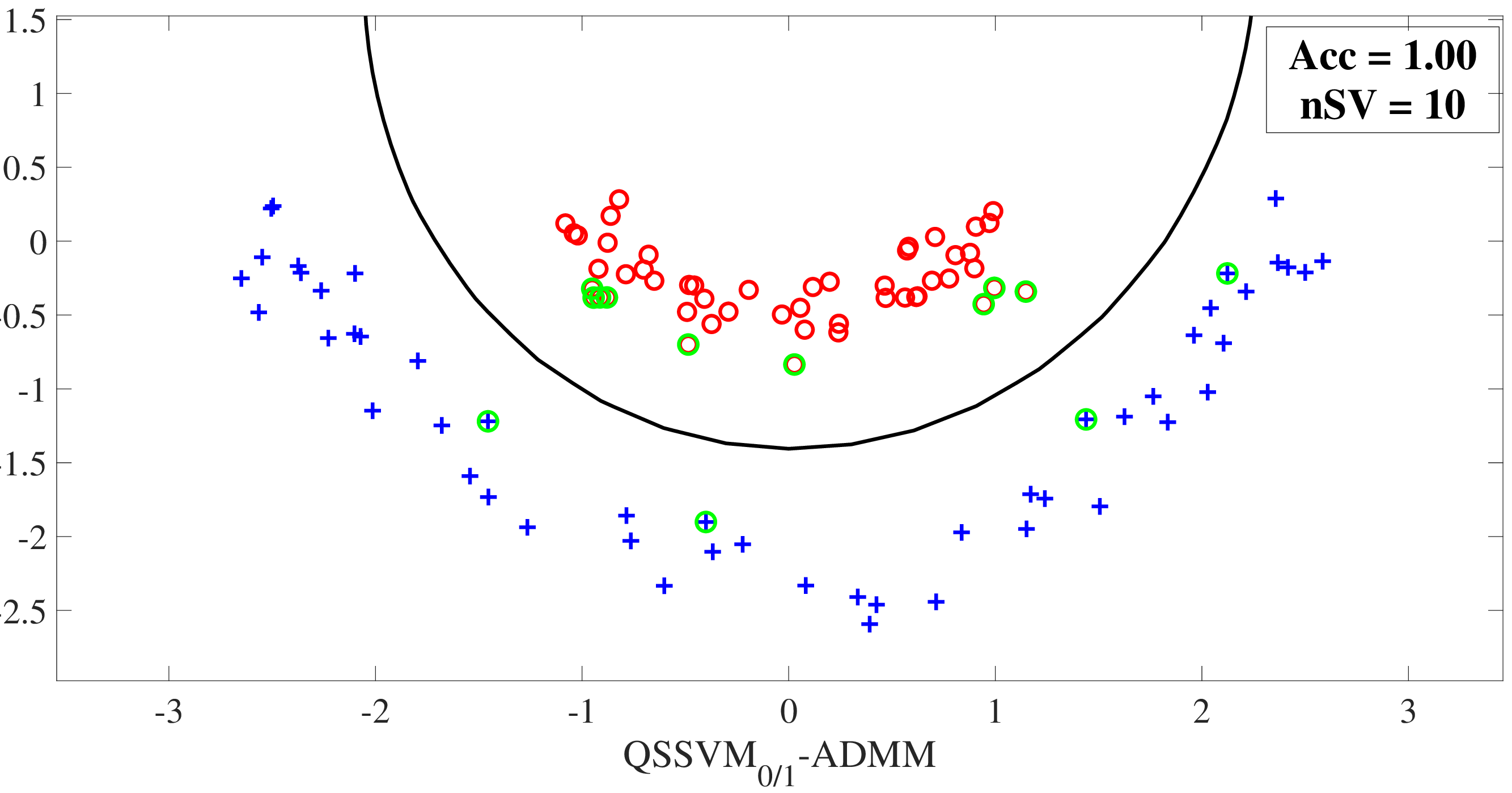}}%}
\hfil
\subfloat[]{\includegraphics[width=1.7in, height= 1in]{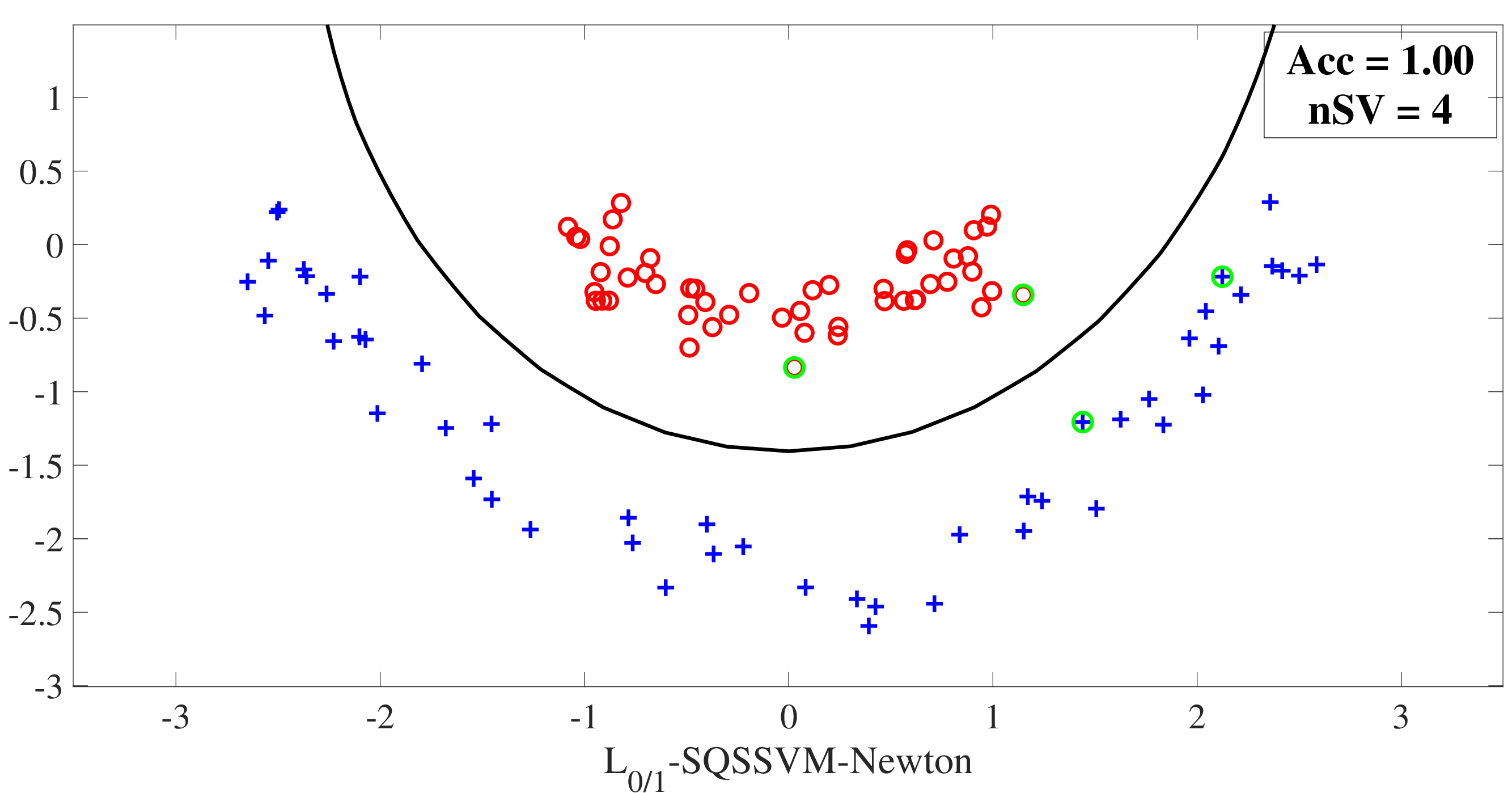}}\\
\subfloat[]{\includegraphics[width=1.7in, height= 1in]{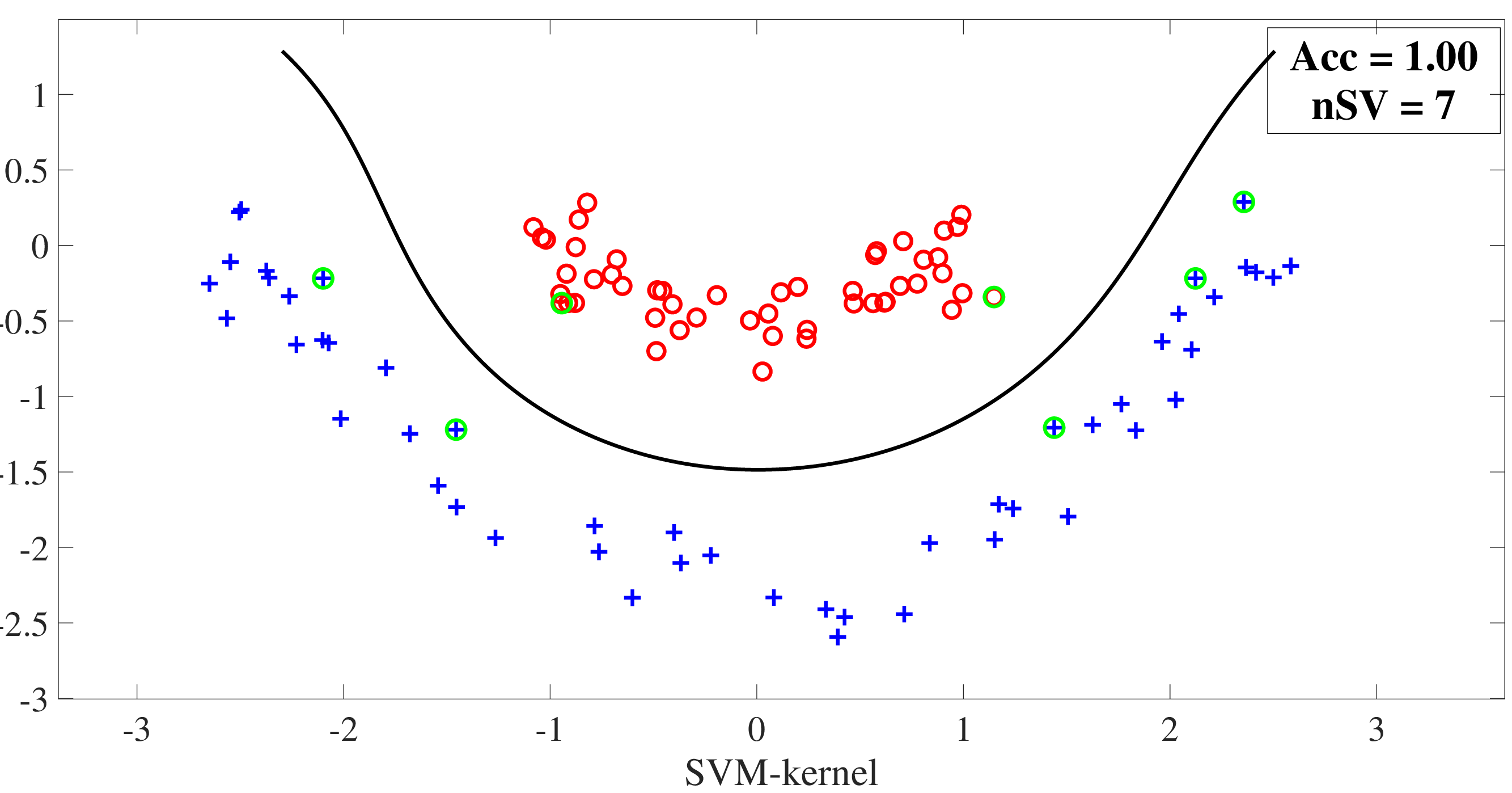}}%
\hfil
\subfloat[]{\includegraphics[width=1.7in, height= 1in]{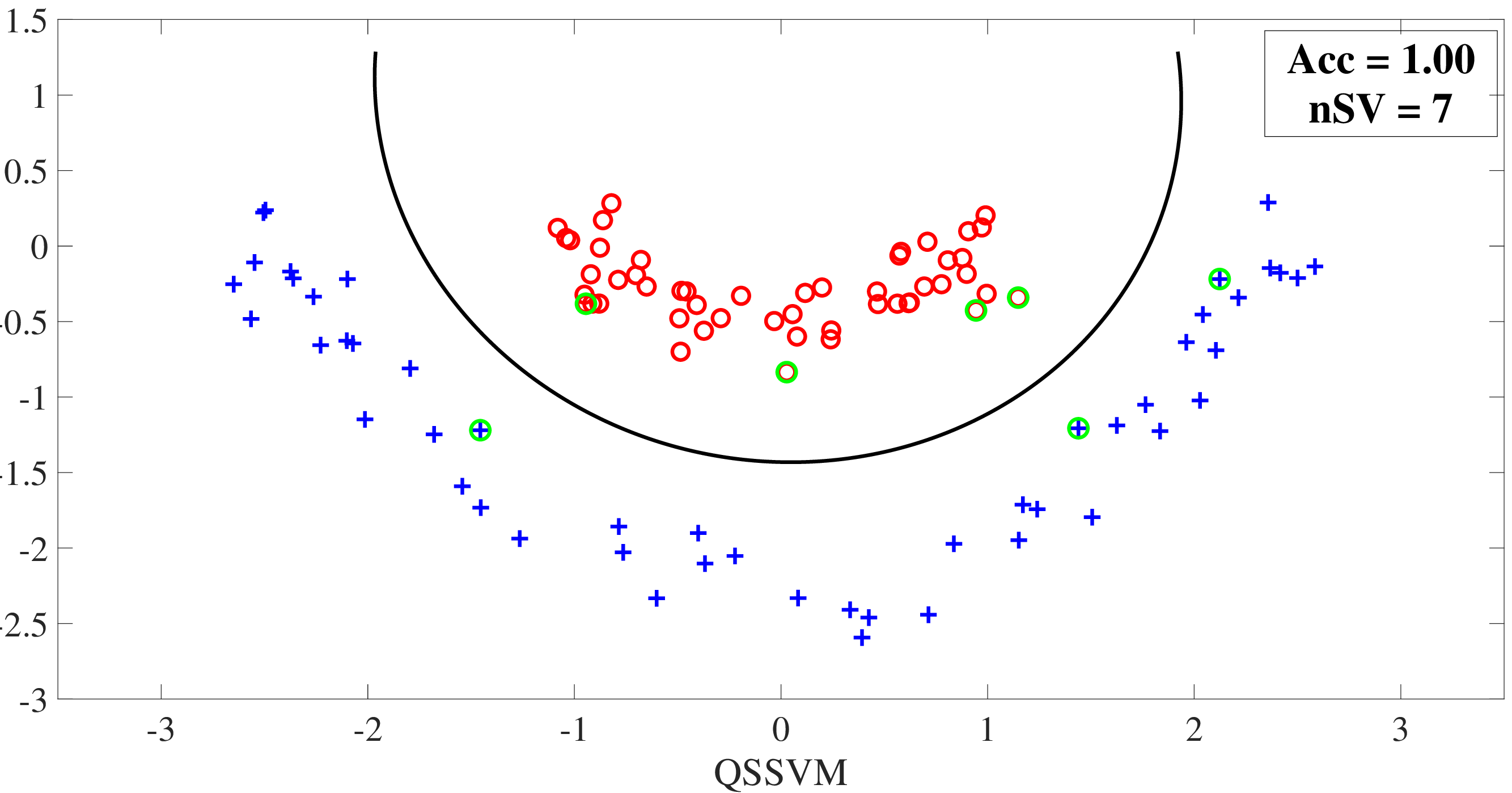}}%
\caption{2d-convex data set  with 100 points.}
\label{Qdata}
\end{figure}

\begin{figure}[H]
\centering
\subfloat[]{\includegraphics[width=1.7in, height= 1in]{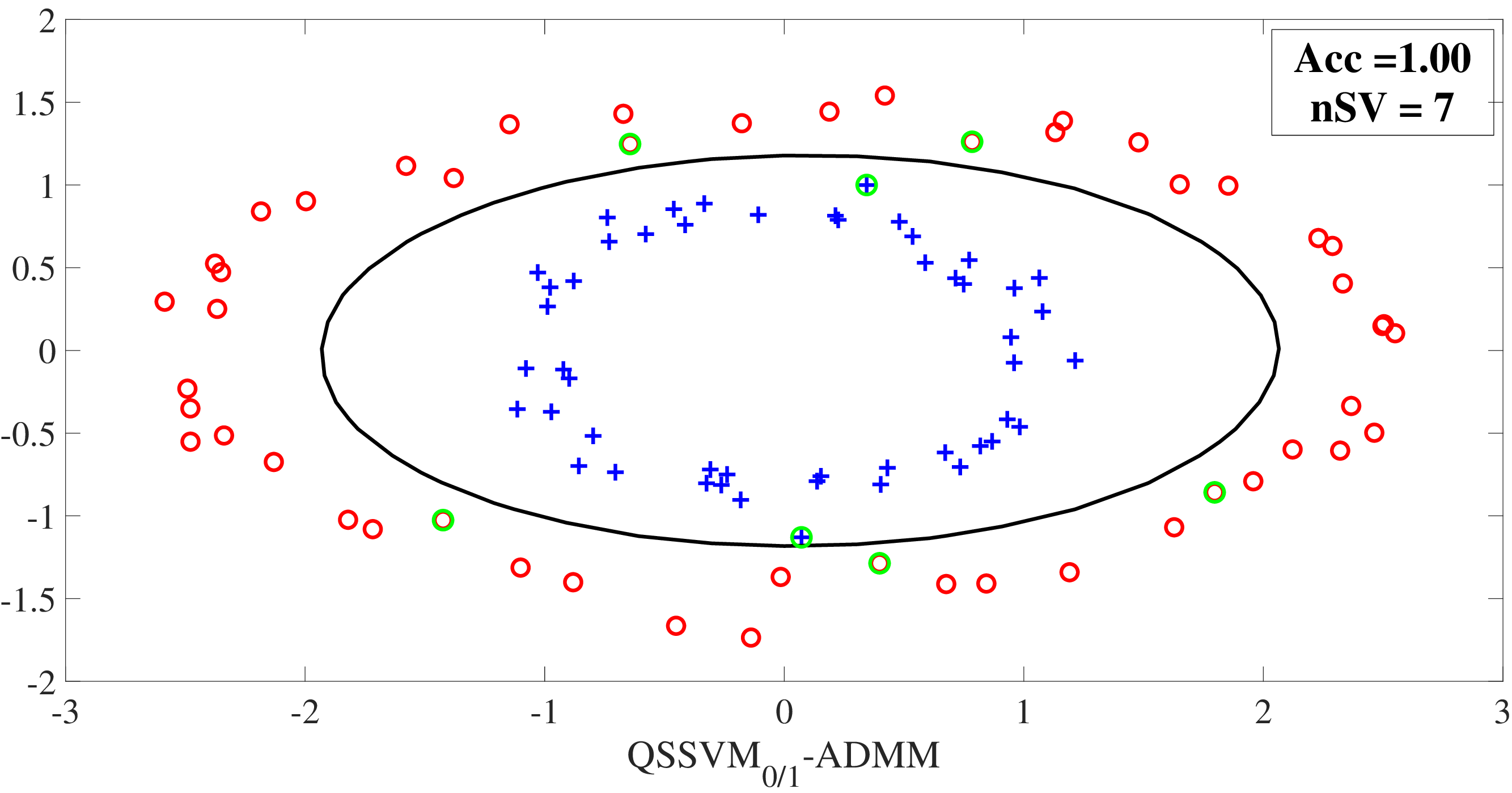}}%}
\hfil
\subfloat[]{\includegraphics[width=1.7in, height= 1in]{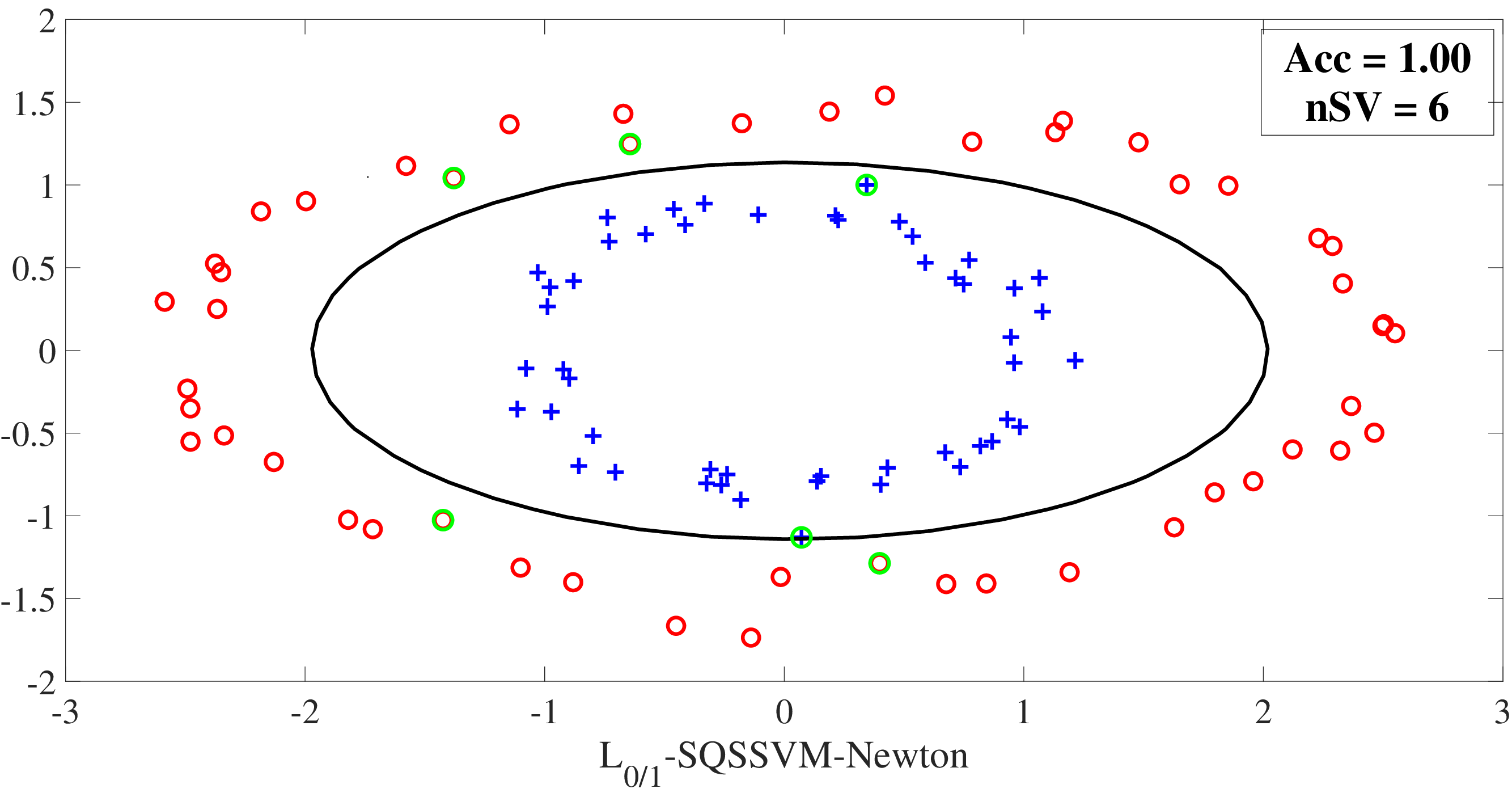}}\\
\subfloat[]{\includegraphics[width=1.7in, height= 1in]{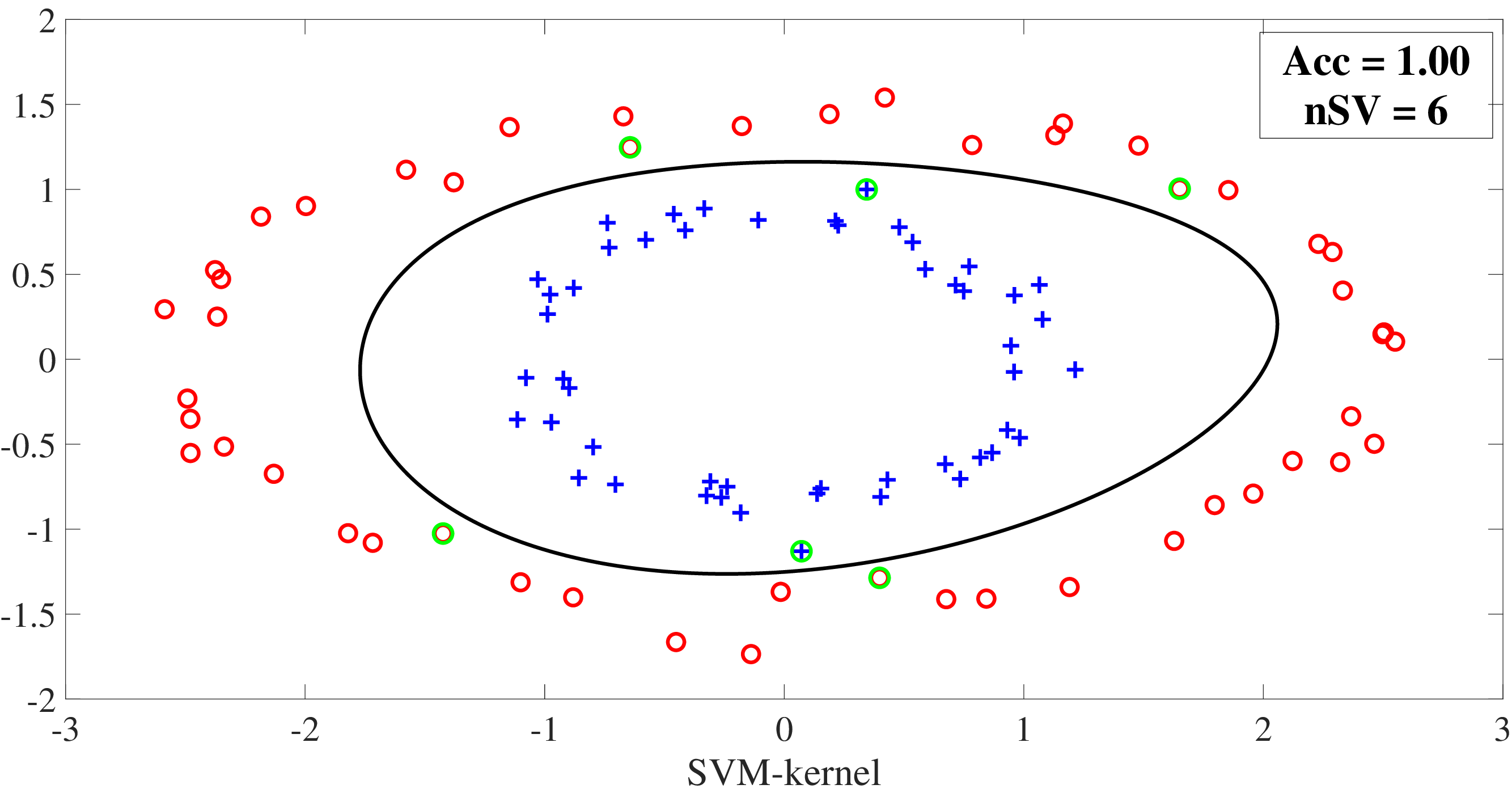}}%
\hfil
\subfloat[]{\includegraphics[width=1.7in, height= 1in]{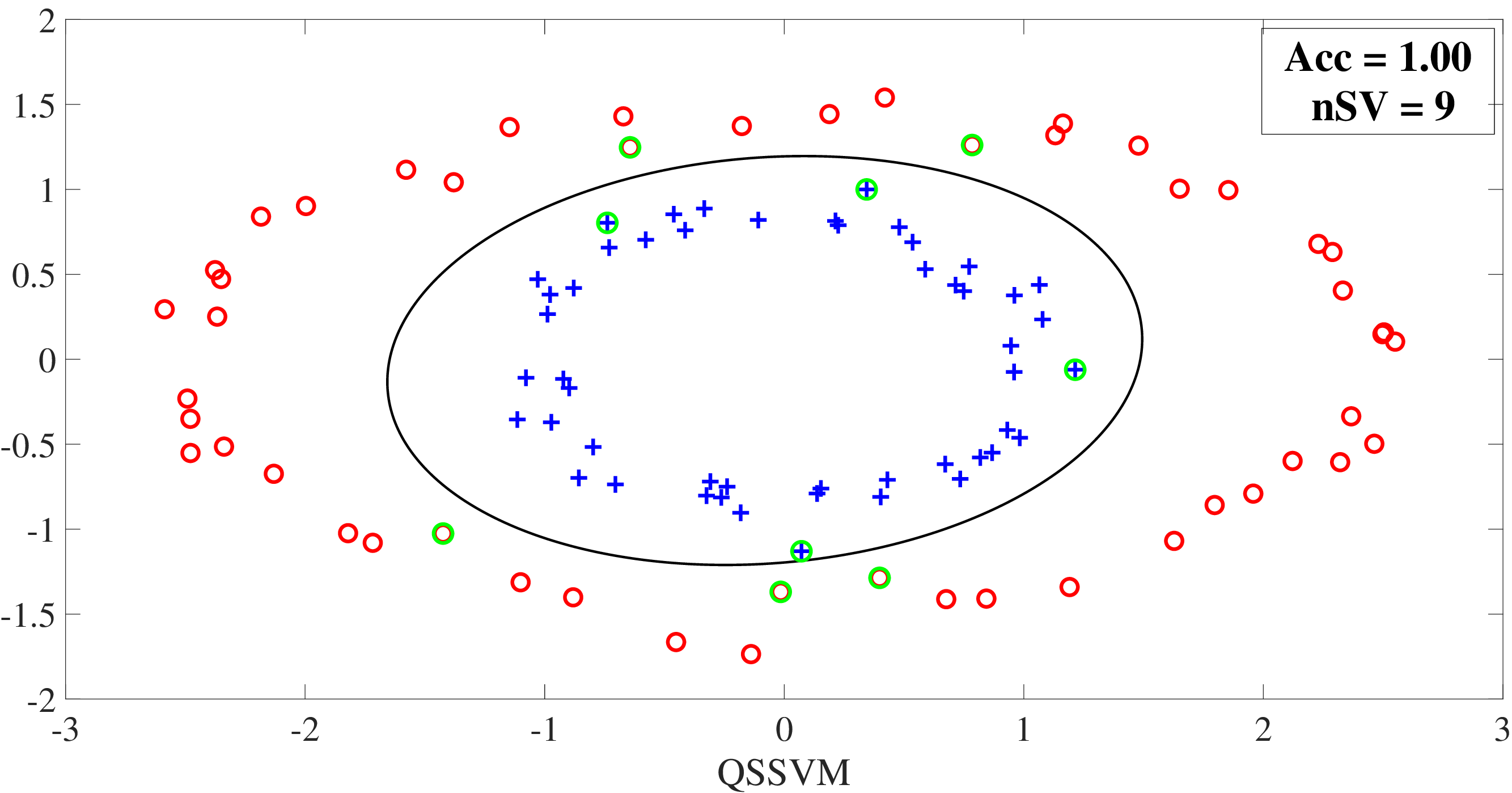}}%
\caption{2d-convex data set  with 100 points.}
\label{2d-convex}
\end{figure}

We first evaluate the performance of $L_{0/1}$-SQSSVM on three artificially generated datasets: a linearly separable data set, a circular (quadratically separable) data set, and a 2d-convex data set. Each data set contains 50 points per class, resulting in 100 points total. The distribution of these datasets is illustrated in the figures. We compare the classification accuracy and CPU time of these methods. The classification results are shown in Figure \ref{lindata} - Figure \ref{2d-convex} and the CPU time is shown in Table \ref{time}  

Figure \ref{lindata} - Figure \ref{2d-convex} depict the flexibility of $L_{0/1}$-SQSSVM in capturing linear and quadratic separating surfaces. As shown in Figure \ref{lindata}, the smooth quadratic hyper-surface obtained by $L_{0/1}$-SQSSVM reduce to the standard SVMs over  linearly separable data sets. In  Figure \ref{Qdata}, \ref{2d-convex} the smooth quadratic hyper-surface of $L_{0/1}$-SQSSVM is the parabola and the circle.  We can see that our method has the lowest number of support vectors (nSV) in the case where the classification accuracy is 1. We compare the minimum,  maximum, and average  of CPU time of the four methods for 50 experiments and  record the values in Table \ref{time}, which shows that our method takes less CPU time compared to the other three methods at the same classification rate. 
\begin{table}[H]
	\caption{The minimum, maximum and average of CPU time.}
	\label{time}
	\centering
	\scriptsize
	\begin{tabular}{ccccccc}
	\hline
	Data sets &   & Min   & Max & Mean  \\
	\hline
	& $L_{0/1}$-SQSSVM-Newton & 0.0024sec &0.0175sec &0.0036sec \\
	Linear& QSSVM$_{0/1}$-ADMM &0.0240sec &0.0680sec &0.0280sec \\
	& SVM with kernel& 0.0021sec &0.2060sec & 0.0102sec \\
	& SQSSVM & 0.0022sec &0.0523sec &0.0038sec \\
	\hline
	& $L_{0/1}$-SQSSVM-Newton & 0.0026sec &0.0183sec &0.0037sec\\
	Circular  & QSSVM$_{0/1}$-ADMM & 0.0078sec &0.0451sec &0.0099sec\\
	& SVM with kernel & 0.0021sec &0.2138sec &0.0104sec \\
	& SQSSVM & 0.0023sec &0.0520sec & 0.0038sec\\
	\hline
	& $L_{0/1}$-SQSSVM-Newton & 0.0014sec & 0.0135sec &0.0022sec\\
	2d-convex& QSSVM$_{0/1}$-ADMM &0.0238sec &0.0668sec & 0.0286sec\\
	& SVM with kernel &0.0021sec &0.2092sec & 0.0102sec\\
	& SQSSVM &  0.0025sec &0.0511sec &0.0041sec\\		
	\hline
	\end{tabular}
\end{table}

To explore the effects of parameters $\lambda$ and $\tau$ on the accuracy of our $L_{0/1}$-SQSSVM, we test the above three data sets. The results are shown in the figure below.
\begin{figure}[H]
\centering
\subfloat[]{\includegraphics[width=1.7in]{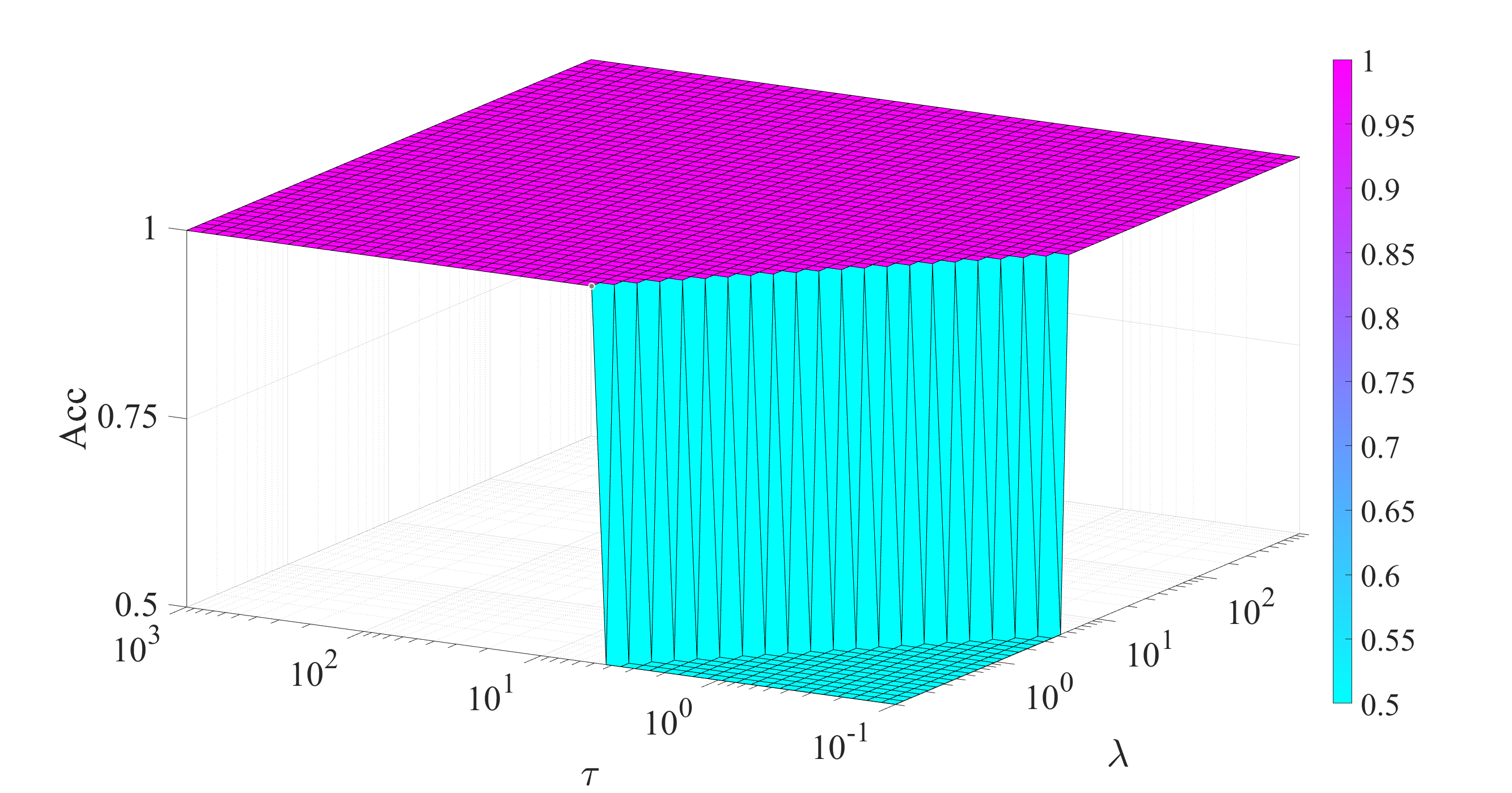}}%}
\hfil
\subfloat[]{\includegraphics[width=1.7in]{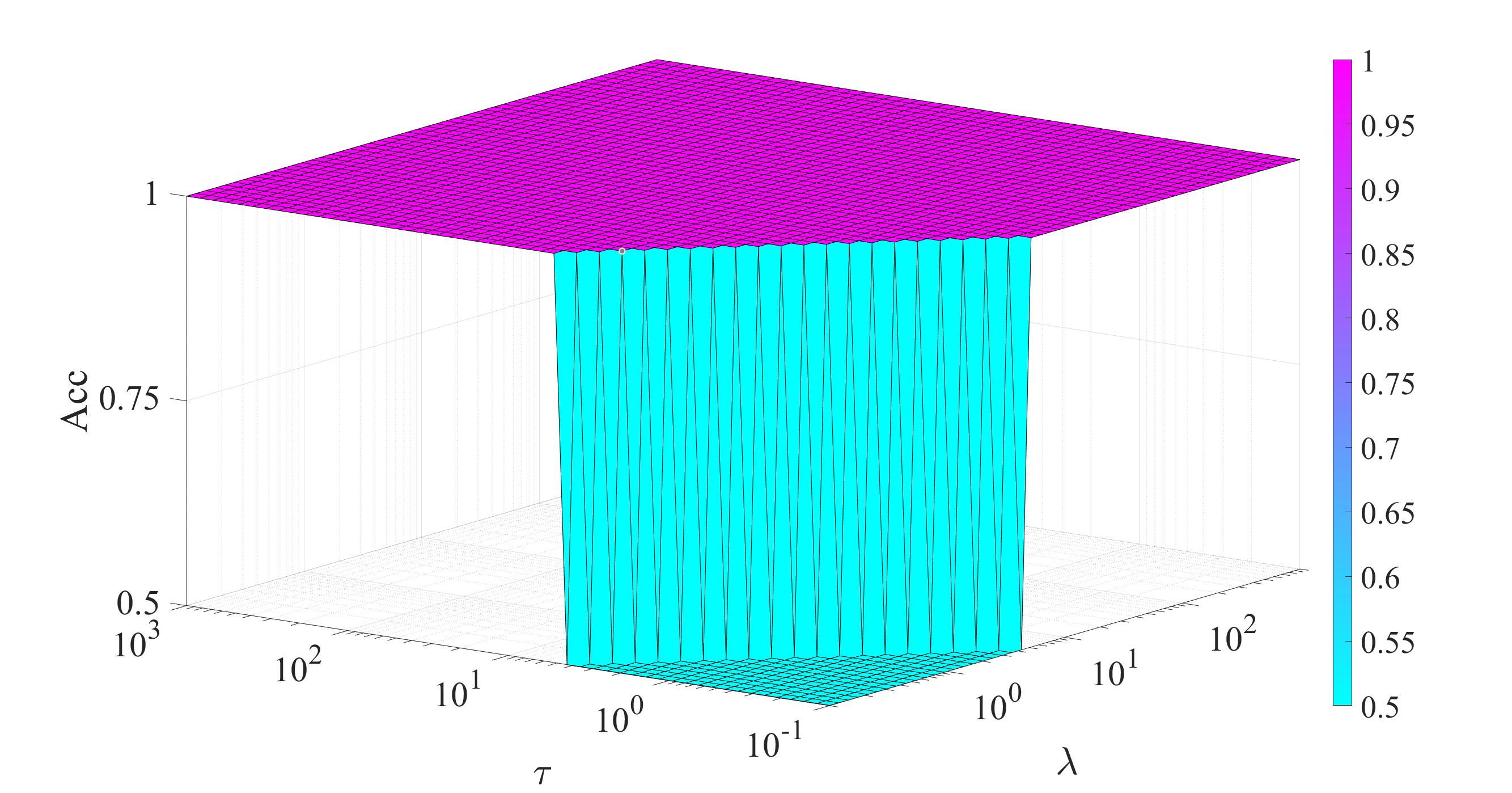}}\\
\subfloat[]{\includegraphics[width=1.7in]{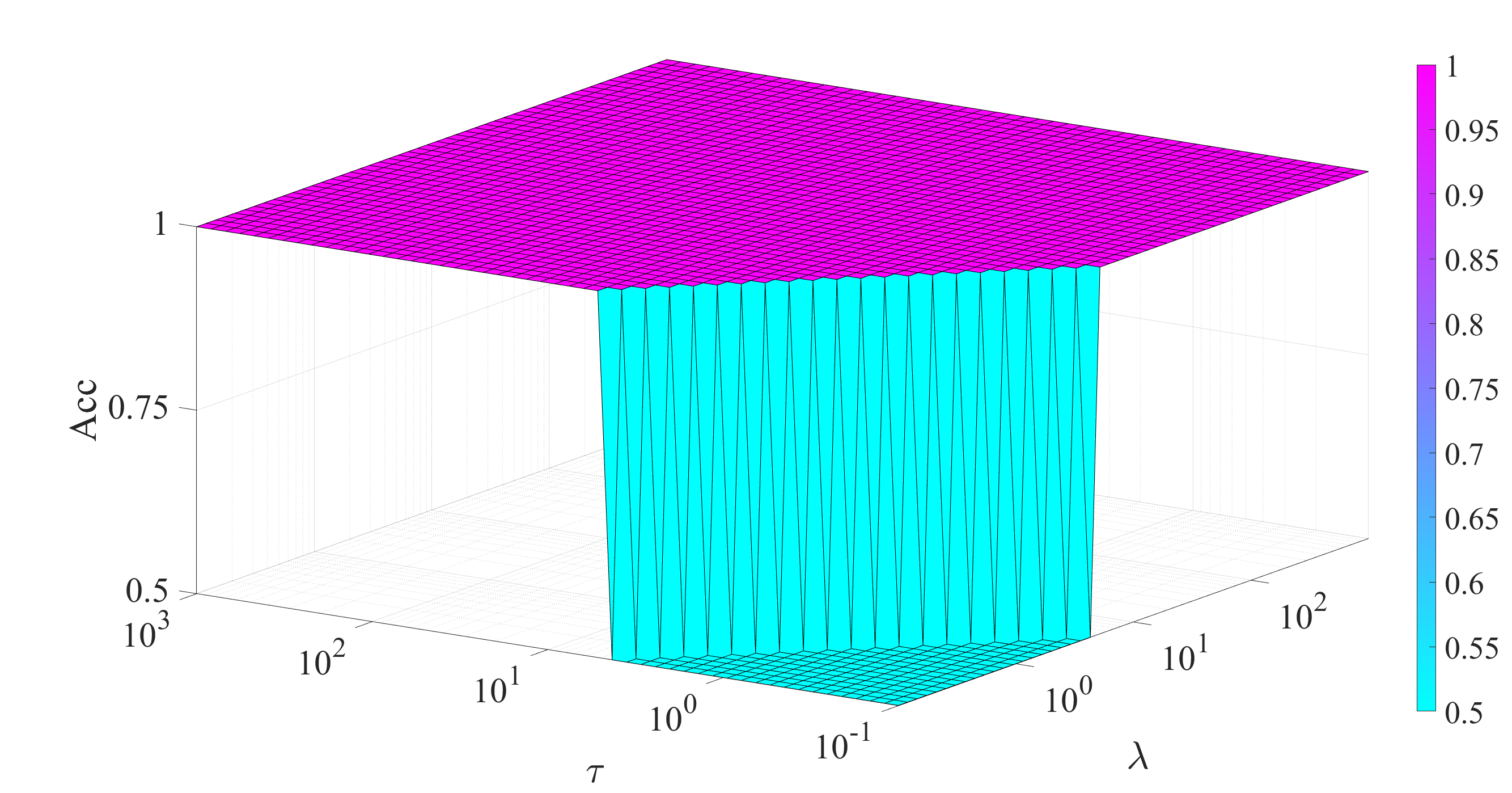}}%
\caption{The accuracy of $L_{0/1}$-QSSVM versus the parameters $\lambda$ and $\tau$ on artificial data sets.}
\label{Qdata2}
\end{figure}
Fig. 4 illustrates that the accuracy is remains 1  when  $\lambda$ and $\tau$ are greater than $10$. These findings suggest that we should take the value of $\lambda$ and $\tau$  to be greater than 1 in the following experiments. 

\subsection{Test with public benchmark data sets}
In this part, a comparative evaluation is conducted on some public benchmark datasets, where $L_{0/1}$-SQSSVM with Newton method is compared with the QSSVM$_{0/1}$ with ADMM, SVM with polynomial kernel and SQSSVM. The information of  datasets are presented in Table \ref{Datasets}.
    \begin{table}[H]
		\caption{Description of  data sets used.} 
		\label{Datasets}
		\centering
		\begin{tabular}{ccccccc}	
		\hline
		Data sets & Features & Classes  & Data size \\ % & Experiment classes & Sample size \\
		\hline
		Iris & 4 & 3 & 150  \\ % & Ver /Vir  &  50/50\\
		\hline
		Wine & 13 & 3 &178 \\ % & Class 1/Class 3 &  59/48\\
		\hline
		Diabetes& 8 & 2 & 768 \\ % &  Yes /No  & 268/500 \\
		\hline
		Balance Scale & 4 & 3 & 625\\ % & B/L & 288/288 \\
		\hline
		Seed & 7 & 3 & 210 \\ % & Kama/Rosa & 70/70 \\
		\hline
		\makecell{Online shoppers\\purchasing  intention}  & 17& 2 & 12330 \\ % &1 /0 & 1908/10422\\
		\hline
		Skin-noskin & 3& 2 &  245057 \\ % & 1/2 & 50859/194198\\
		\hline
		Htru2 & 8& 2 & 17898 \\ % & 1/0 & 1639/16259\\
		\hline		
		\end{tabular}
	\end{table}

We can see that there are more than two classes for each data. In later  experiments, we choose only two classes for binary classification. For each model, there are different percentage of training data and  the experiments are repeated for 50 times. The minimum,  maximum,  mean and standard deviation  of accuracy scores  and the average CPU time among these 50 experiments are recorded in the table, where $"-"$ denotes the results are not obtained if a solver takes too much time or requires a large memory that is out of the capacity of our desktop.  The following Table \ref{Iris}-\ref{Seed}  show that the proposed $L_{0/1}$-SQSSVM  with Newton method is  more efficient than the other methods for benchmark data sets.\\
(a) The mean and variance of accuracy scores obtained by $L_{0/1}$-SQSSVM are  better than those of other methods over most of  the tested benchmark data sets for different training rate. This indicates that the proposed method has the classification stability on real data sets.\\
(b) The CPU times of $L_{0/1}$-SQSSVM on the tested data sets are acceptable and less than QSSVM$_{0/1}$ with ADMM in most cases. \\
(c) For large datasets, while other methods may achieve higher  mean and variance of accuracy scores, our method achieves significantly superior time performance.
\begin{table}[H]
	\caption{Results for the  Iris data set}
	\label{Iris}
	\centering
	\tiny
	\begin{tabular}{ccccccc}
	\hline
	Training rate &Model & Min(k\%) & Max(k\%)  & Mean(k\%) & Var & CPU time \\
	\hline
	& $L_{0/1}$-SQSSVM  &46.25 & 97.50 & 89.93 &  0.094 & 0.004sec\\		
	& QSSVM$_{0/1}$  &36.25 & 93.75 & 53.98 &  0.098 & 0.032sec\\
	20\%  & SQSSVM & 75.00 &  97.50 & 91.88 &   0.045 & 0.003sec  \\ 
	& SVM & 73.75 & 96.25 & 88.15  &  0.045 & 0.003sec\\
	\hline
	&  $L_{0/1}$-SQSSVM  &48.33 & 100.00 & 91.97 &  0.073 & 0.010sec \\
	& QSSVM$_{0/1}$  &36.67 & 88.33 & 51.37 & 0.083 & 0.030sec\\
	40\% & SQSSVM & 83.33 &  98.33 &92.47 &   0.039 & 0.003sec\\
	& SVM& 83.33 &  95.00 & 90.36 &    0.039 & 0.002sec\\
	\hline
	& $L_{0/1}$-SQSSVM & 85.00 & 100.00 & 92.75 &  0.041 & 0.019sec\\
	& QSSVM$_{0/1}$  &  32.50 & 97.50 & 52.55 & 0.133 & 0.033sec\\
	60\%& SQSSVM &  80.00 &100.00 & 93.95 & 0.037 & 0.004sec \\
	& SVM &   80.00 & 97.50 & 90.85 & 0.053 & 0.003sec \\
	\hline
	& $L_{0/1}$-SQSSVM &84.21 & 100.00 &93.26 &  0.048 &  0.030sec\\
	& QSSVM$_{0/1}$  &26.32 & 100.00 & 56.00 & 0.164 & 0.034sec\\
	80\% & SQSSVM& 73.68 & 100.00 &96.21 &   0.061 & 0.006sec \\
	& SVM&  68.42 & 100.00& 90.21 &  0.036 & 0.003sec \\			
	\hline
	\end{tabular}
\end{table}

	\begin{table}[H]
		\caption{Results for the  Wine data set}
		\label{Wine}
		\centering
		\tiny
		\begin{tabular}{ccccccc}
		\hline
		Training rate &Model & Min(k\%) & Max(k\%)  & Mean(k\%) & Var & CPU time \\
		\hline
		& $L_{0/1}$-SQSSVM  &89.41 & 100.00 & 97.86 & 0.027 & 0.016sec\\
		20\%& QSSVM$_{0/1}$  &0.00 &  67.06 & 11.91 & 0.215 & 0.080sec\\
		& SQSSVM &85.88 & 94.12 & 90.00 & 0.017 & 0.003sec \\
		& SVM & 84.71 &100.00 &97.81 & 0.017 & 0.002sec\\
		\hline
		&  $L_{0/1}$-SQSSVM  & 98.44 & 100.00 & 99.94& 0.003 & 0.026sec\\
		& QSSVM$_{0/1}$  &0.00 & 87.50 & 41.81 & 0.204 &  0.157sec\\
		40\% & SQSSVM & 81.25 & 96.88 & 89.94 & 0.030 & 0.003sec\\
		& SVM&96.88 &  100.00 & 99.63 & 0.030 & 0.002sec\\
		\hline
		& $L_{0/1}$-SQSSVM &100.00 & 100.00 &  100.00  & 0.000 & 0.036sec \\
		& QSSVM$_{0/1}$  & 33.33 & 95.24 & 50.67 & 0.116 & 0.193sec\\
		60\%& SQSSVM &  76.19 & 97.62 & 89.71 & 0.043 & 0.003sec \\
		& SVM &   97.62 & 100.00 & 99.86 & 0.043 & 0.002sec\\
		\hline			
		& $L_{0/1}$-SQSSVM & 100.0 & 100.00 & 100.0 & 0.000 & 0.045sec\\
		& QSSVM$_{0/1}$  & 28.57 & 100.00 & 52.19 & 0.140 &  0.201sec\\
		80\% & SQSSVM& 71.43 &  100.00 &  90.48 &   0.062 & 0.007sec \\
		& SVM&  95.24 &  100.00 & 99.90 &   0.062 & 0.003sec \\			
		\hline
		\end{tabular}
	\end{table}

	\begin{table}[H]
		\caption{Results for the  Diabetes data set}
		\label{Diabetes}
		\centering
		\tiny
		\begin{tabular}{ccccccc}
		\hline
		Training rate &Model & Min(k\%) & Max(k\%)  & Mean(k\%) & Var & CPU time \\
		\hline
		& $L_{0/1}$-SQSSVM  & 73.81 & 78.09 & 75.89 &  0.009 & 0.064sec\\ 
		20\%  & QSSVM$_{0/1}$  &47.94 & 72.65 & 66.57 & 0.047 & 0.087sec\\
		& SQSSVM & 72.16 & 78.75 &75.54 & 0.011 & 0.013sec \\
		& SVM& 56.67 & 70.51 & 65.04 & 0.011 & 0.006sec\\
		\hline
		&  $L_{0/1}$-SQSSVM  &  73.85 & 80.22 & 76.83 & 0.013 & 0.116sec\\
		& QSSVM$_{0/1}$  &62.42 & 75.82 & 68.23 & 0.033 & 0.119sec\\
		40\% & SQSSVM & 72.31 &79.78 & 76.12 & 0.013 & 0.008sec\\
		& SVM& 61.10 & 72.31 & 65.96 & 0.023 & 0.033sec\\
		\hline
		& $L_{0/1}$-SQSSVM &72.23 & 83.50 & 76.94 & 0.020 & 0.168sec\\
		& QSSVM$_{0/1}$  &62.05 &76.24 & 68.07 & 0.033 &  0.160sec\\
		60\% & SQSSVM & 72.23 &  83.50 & 76.09 & 0.023 & 0.014sec\\
		& SVM & 61.39 & 73.60 & 66.83 & 0.023 & 0.440sec\\
		\hline
		& $L_{0/1}$-SQSSVM & 70.86 & 85.43 & 76.93 & 0.031 & 0.225sec\\
		& QSSVM$_{0/1}$  &59.60 & 82.12 & 67.83 & 0.044 & 0.189sec\\
		80\% & SQSSVM & 68.87 & 82.12 & 76.23 & 0.028 & 0.020sec\\
		& SVM &56.95 & 75.50 & 67.25 & 0.028 &  2.104sec\\
		\hline
		\end{tabular}
		\end{table}
		
	\begin{table}[H]
		\caption{Results for the Balance Scale data set}
		\label{Balance}
		\centering
		\tiny
		\begin{tabular}{ccccccc}
		\hline
		Training rate &Model & Min(k\%) & Max(k\%)  & Mean(k\%) & Var & CPU time \\
		\hline
		& $L_{0/1}$-SQSSVM  &  92.61 & 99.13 & 96.65 & 0.013 & 0.035sec\\ 
		& QSSVM$_{0/1}$  &60.43 & 96.52 & 82.49 &  0.068 & 0.035sec\\
		20\%  & SQSSVM & 92.83 & 96.52 & 94.74 & 0.008 & 0.003sec\\
		& SVM& 90.22 &  98.48 & 95.06 & 0.008 & 0.003sec\\
		\hline
		&  $L_{0/1}$-SQSSVM  & 94.49 &100.00 & 97.33 & 0.012 & 0.076sec\\
		& QSSVM$_{0/1}$  &72.17 & 95.65 & 85.99 &  0.048 & 0.044sec \\
		40\% & SQSSVM & 92.75 & 97.97 & 94.53 &  0.009 & 0.003sec \\
		& SVM& 94.78 & 98.84 & 97.33 & 0.009 & 0.004sec\\
		\hline
		& $L_{0/1}$-SQSSVM & 95.22 & 100.00 & 97.95 & 0.011 & 0.120sec\\
		& QSSVM$_{0/1}$  &76.96 & 97.39 & 87.44 & 0.051 & 0.052sec\\
		60\% & SQSSVM & 90.87 & 97.39 & 94.57 &  0.014 & 0.004sec\\
		& SVM & 96.09 &  99.57 &  97.83 &  0.014 & 0.006sec\\
		\hline
		& $L_{0/1}$-SQSSVM & 93.91 & 100.00 & 98.14 & 0.016 & 0.154sec\\
		& QSSVM$_{0/1}$  &76.52 & 98.26 & 89.15 & 0.047 & 0.059sec\\
		80\% & SQSSVM &89.57 & 98.26 & 94.19 &  0.019 & 0.004sec\\
		& SVM & 94.78 & 99.13 & 97.91 &  0.019 & 0.010sec\\
		\hline
		\end{tabular}
	\end{table}

	\begin{table}[H]
		\caption{Results for the Seed data set}
		\label{Seed}
		\centering
		\tiny
    	\begin{tabular}{ccccccc}
		\hline
		Training rate &Model & Min(k\%) & Max(k\%)  & Mean(k\%) & Var & CPU time \\
		\hline
		& $L_{0/1}$-SQSSVM  & 90.18 &  100.00 & 97.38 & 0.018 & 0.002sec\\ 
		& QSSVM$_{0/1}$  & 33.93 & 84.82 & 51.27 & 0.071 & 0.067sec\\ 
		20\%  & SQSSVM & 90.18 & 96.43 & 93.07 & 0.014 & 0.020sec\\
		& SVM& 81.25 & 96.43& 90.18 &  0.014 & 0.004sec\\
		\hline
		&  $L_{0/1}$-SQSSVM  & 92.86 & 100.00 & 97.67 & 0.014 & 0.005sec\\
		&QSSVM$_{0/1}$ &38.10 & 90.48 & 52.79 & 0.095 & 0.073sec\\
		40\% & SQSSVM &  86.90 & 98.81 & 93.14 & 0.024 & 0.003sec \\
		& SVM &  81.25 & 96.43 & 90.20 & 0.024 & 0.003sec\\
		\hline
		& $L_{0/1}$-SQSSVM & 94.64 &100.00 & 98.04 & 0.016 & 0.008sec\\
		&QSSVM$_{0/1}$ & 37.50 & 94.64 & 53.86 & 0.104 & 0.077sec\\
		60\% & SQSSVM & 87.50 & 100.00 & 93.32 & 0.027 & 0.003sec\\
		& SVM &  85.71 & 100.00 & 93.90 & 0.027 &  0.003sec\\
		\hline
		& $L_{0/1}$-SQSSVM &92.59 &100.00 &98.30 &  0.024 & 0.008sec\\
		& QSSVM$_{0/1}$  & 33.33 & 96.30&  53.85  &  0.120 & 0.079sec\\
		80\% & SQSSVM &77.78 & 100.00 & 93.70 & 0.048 & 0.003sec\\
		& SVM &  85.19 & 100.00 & 95.04 & 0.048 & 0.003sec\\
		\hline
		\end{tabular}
	\end{table}

	\begin{table}[H]
		\caption{Results for Online shoppers purchasing intention data set}
		\centering
		\tiny
		\begin{tabular}{ccccccc}
		\hline
		Training rate &Model & Min(k\%) & Max(k\%)  & Mean(k\%) & Var & CPU time \\
		\hline
		& $L_{0/1}$-SQSSVM  & 87.62 & 88.69 & 88.33 & 0.003 & 0.768sec\\ 
		& QSSVM$_{0/1}$&19.70 &85.29 & 36.90 & 0.187 & 6.640sec\\ 
		20\%  & SQSSVM &87.84 & 88.98 & 88.48 & 0.004 & 2.212sec\\
		& SVM& 80.46 & 82.65 & 81.86 & 0.008 & 2.721sec\\
		\hline
		&  $L_{0/1}$-SQSSVM  &88.05 & 88.69 & 88.37 & 0.002 & 1.613sec\\ 
		&QSSVM$_{0/1}$ &27.26 & 42.74 & 36.90 & 0.047 & 15.308sec\\
		40\% & SQSSVM &88.06 &89.29 & 88.65 & 0.004 & 3.971sec\\ 
		& SVM & 84.10 &84.82 & 84.46 & 0.003 & 45.667sec\\
		\hline
		& $L_{0/1}$-SQSSVM & 88.30 & 88.99 & 88.55 & 0.002 & 2.405sec\\
		&QSSVM$_{0/1}$ &25.51 &50.91 & 36.06& 0.075 & 21.024sec\\
		60\% & SQSSVM & 46.63 & 89.25 &80.94 & 0.155 & 6.158sec\\
		& SVM &85.28 & 86.23 & 85.62 & 0.004 & 85.380sec\\
		\hline
		& $L_{0/1}$-SQSSVM & 88.28 & 89.41 &88.64 & 0.004 & 3.255sec\\ 
		& QSSVM$_{0/1}$  &26.86 & 71.72 & 40.30 & 0.165 & 25.522sec\\
		80\% & SQSSVM &31.93 &88.97 &78.15 & 0.188 & 8.030sec\\
		& SVM & 84.02 & 85.76 & 84.63 & 0.006 & 90.782sec\\
		\hline
		\end{tabular}
		\label{Seed} 
	\end{table} 

	\begin{table}[H]
		\caption{Results for the Skin-noskin data set}
		\label{Seed}
		\centering
		\tiny
		\begin{tabular}{ccccccc}
		\hline
		Training rate &Model & Min(k\%) & Max(k\%)  & Mean(k\%) & Var & CPU time \\
		\hline
		& $L_{0/1}$-SQSSVM  & 95.81 &95.98 & 95.85 & 0.001 & 1.024sec\\ 
		& QSSVM$_{0/1}$  &58.57 &  64.03 & 62.82 &  0.021 & 1.935sec\\ 
		20\%  & SQSSVM &- &- & - & - & -\\
		& SVM &99.84 & 99.88 & 99.86 & 0.000 & 19.948sec\\
		\hline
		&  $L_{0/1}$-SQSSVM  & 95.76 & 95.96 & 95.84 & 0.001 & 2.015sec\\
		&QSSVM$_{0/1}$ &63.08 & 64.03 & 63.71 & 0.003 & 3.849sec\\
		40\% & SQSSVM&- &- & - & - & -\\
		& SVM &99.85 & 99.87 & 99.86  & 0.000 & 53.003sec\\
		\hline
		& $L_{0/1}$-SQSSVM & 95.65 & 95.89 &95.81 & 0.001 & 3.054sec\\
		&QSSVM$_{0/1}$ & 63.63 & 64.10 & 63.90 & 0.002  &12.609sec\\
		60\% & SQSSVM &- &- & - & - & -\\
		& SVM &99.87 & 99.88 & 99.88 & 0.000 & 47.684sec\\
		\hline
		& $L_{0/1}$-SQSSVM &95.85 & 95.91 &95.88 & 0.000 & 3.733sec\\
		& QSSVM$_{0/1}$  &63.61 & 64.28 & 63.93 & 0.002 & 16.614sec\\
		80\% & SQSSVM&- &- & - & - & -\\
		& SVM &99.88 & 99.91 & 99.89 &  0.000 & 180.111sec\\
		\hline
		\end{tabular}
	\end{table} 

	\begin{table}[H] 
		\caption{Results for the Htru2 data set}
		\label{Seed}
		\centering
		\tiny
		\begin{tabular}{ccccccc}
		\hline
		Training rate &Model & Min(k\%) & Max(k\%)  & Mean(k\%) & Var & CPU time \\
		\hline
		& $L_{0/1}$-SQSSVM  &97.26 & 97.45 & 97.38 &0.001 & 0.166sec\\ 
		& QSSVM$_{0/1}$  &10.12 & 95.02 & 76.09  &  0.304 & 0.626sec\\ 
		20\%  & SQSSVM &97.10 & 97.50 & 97.30 & 0.001 &  0.112sec\\
		& SVM&97.21 & 97.62 & 97.51 & 0.002 &  12.501sec\\
		\hline
		&  $L_{0/1}$-SQSSVM  &97.37 & 97.72 & 97.48 & 0.001 & 0.350sec\\
		&QSSVM$_{0/1}$ & 84.78 & 94.62 & 91.26 & 0.026 & 1.100sec\\
		40\% & SQSSVM &97.13 &97.74 & 97.39 & 0.001 & 0.796sec\\
		& SVM &97.73 & 98.01 & 97.84 & 0.001 & 60.277sec\\
		\hline
		& $L_{0/1}$-SQSSVM &97.18 & 97.49 & 97.29 & 0.001 & 0.522sec\\
		&QSSVM$_{0/1}$ &91.31 & 93.73 &  92.50 &  0.009 & 4.523sec\\
		60\% & SQSSVM & 97.09 & 97.82 &97.44  &  0.001 & 2.619sec\\
		& SVM & 97.58 & 97.96 & 97.82 & 0.001 & 98.994sec\\
		\hline
		& $L_{0/1}$-SQSSVM & 96.95 & 97.76 & 97.27 &  0.002 & 0.694sec\\
		& QSSVM$_{0/1}$  &91.98 & 93.77 & 92.83 & 0.005 & 4.169sec\\
		80\% & SQSSVM &96.90 & 97.90 & 97.47 &0.002 & 4.517sec\\
		& SVM &97.54 & 98.07 & 97.87 & 0.002 & 123.251sec\\
		\hline
	\end{tabular}
	\end{table} 
	
\section{Conclusion}
In this paper, we proposed a sparsity optimization problem \eqref{Prim}, which is the prototype of the kernel-free soft-margin  quadratic surface support vector machine model with 0-1 loss ($L_{0/1}$-SQSSVM). We demonstrate the effectiveness of our method through both theoretical analysis and numerical experiments. Firstly, we derived the first-order and second-order optimality conditions for the problem by utilizing the concept of subderivatives, and then we developed the Newton method using the P-stationary equation and  established its local quadratic convergence under some reasonable assumptions. Next, we conducted extensive numerical experiments to show the practical efficiency of the proposed method. 

The experimental results demonstrate  that our $L_{0/1}$-SQSSVM  with Newton method is  not only suitable for quadratically separable datasets but also performs robustly on linearly separable datasets. For real datasets, our method also performs favorably compared to other methods.

\section*{Acknowledgments}


\bibliographystyle{IEEEtran}
\bibliography{refer}% common bib file
	
	{\appendix
		\section*{Proof of  Proposition \ref{limsubdiff}}
		From Lemma \ref{lemma2.2}, we have that $g$ is regular at a given $ F(X^*) \in \R^n$ and 
		\begin{align}\label{regular}
			\partial^pg( F(X^*))&=\hat{\partial}g( F(X^*))=\partial g( F(X^*))=\partial^{\infty}g( F(X^*)) \nonumber\\
			&=\{v \in  \mathbb{R}^n ~|~v_i \ge 0,~i \in \mathcal{T}_*  ;~v_i =0,~i \notin\mathcal{T}_*\}.
		\end{align}
		For every $v \in \partial^{\infty} g(F(X^*))=\{v \in \mathbb{R}^n~|~v_i \geq 0, i \in{\mathcal{T}_*};
		~v_i =0, i \not\in {\mathcal{T}_*}\}$ with $A^Tv =0$, the Assumption \ref{Ass2.1} implies that $v=0$.
		It follows from  \cite[Proposition 10.6]{RW} and $g$ is regular at $F(X^*)$ that  $\varphi$ is regular at $X^*$ and
		\begin{align*}
			\partial \varphi(X^*) =A^T \partial \|F(X^*)_+\|_0,
		\end{align*}
		\begin{align*}
			{\rm d}\varphi(X^* )(H)  =	{\rm d} g(F(X^*))(AH).
		\end{align*}
		From Lemma \ref{lemma2.2}, we have 
		\begin{align*}
			{\rm d} g(F(X^*))(AH) = 0, 
		\end{align*}	
		for $AH \in {\rm dom}~{\rm d}g(X^* ) = \{AH\in  \mathbb{R}^n~|~\langle A_i, H\rangle\leq 0,~i \in \mathcal{T}_*;~\langle A_i, H\rangle\in \mathbb{R},~i \notin \mathcal{T}_*\}.$  Therefore, 
		${\rm dom}~{\rm d}\varphi(X^* ) = \{H\in  \mathbb{R}^{m\times p}~|~\langle A_i,H\rangle \leq 0,~i \in \mathcal{T}_*;~\langle A_i, H\rangle \in \mathbb{R},~i \notin \mathcal{T}_*\}$ and 
		\begin{align*}
			{\rm d}\varphi(X^*)(H) 
			=\left\{\begin{array}{clc}
				0, &H \in {\rm dom}~{\rm d}\varphi(X^*),\\ 
				+\infty, &{\rm otherwise}.
			\end{array}
			\right.
		\end{align*}
		
		\section*{Proof of Proposition \ref{prop-critica-varphi}}
		By the definition of the critical cone  in \eqref{critical}, we have 
		\begin{equation*}
			K_{\varphi}(X^*,V) = \{ H \in \mathbb{R}^{m \times p} ~|~ \langle V,H \rangle ={\rm d}\varphi(X^*)(H)\}.
		\end{equation*}
		For any $V \in \partial \varphi(X^*)$, there exists $z^* \in\Lambda(X^*, V)$ such that $ V = A^Tz^*$. Therefore,  we have
		\small{\begin{align*}
				{\rm d}\varphi(X^*)(H)&=\langle V,H \rangle = \langle A^Tz^*, H \rangle=\langle z^*, AH \rangle\\
				&=\underset{i \in \mathcal{T}_*}{\sum}z_i^*\langle  A_i,H \rangle +\underset{i \notin \mathcal{T}_*}{\sum}z_i^*\langle  A_i,H \rangle\\
				&=\underset{i \in \mathcal{T}_*}{\sum}z_i^*\langle  A_i,H\rangle .
		\end{align*}}
		The last equality due to $z_i^* =0,~ i \in \bar{\mathcal{T}}_*$. Further,
		\begin{equation*}
			K_{\varphi}(X^*,V) = \{ H \in \mathbb{R}^{m \times p} ~|~{\rm d}\varphi(X^*)(H)= \underset{i \in \mathcal{T}_*}{\sum}z_i^*\langle  A_i,H\rangle\}.
		\end{equation*}
		Denote by $S = \left\{D  \in \R^{m\times p}~|~\begin{aligned}
			&z^*_i \langle  A_i,D\rangle=0,~\langle  A_i,D\rangle \leq 0,\\
			&z^*_i\geq 0,~~i\in \mathcal{T}_*
		\end{aligned}\right\}$. On the one hand, for any $H \in K_{\varphi}(X^*,V)$, we have that
		\begin{align*}
			{\rm d}\varphi(X^*)(H) = \underset{i \in \mathcal{T}_*}{\sum}z_i^*\langle  A_i,H\rangle
		\end{align*}
		is finite. By Proposition \ref{limsubdiff}, $H \in {\rm dom~ d}\varphi(X^*)=\{H~|~\langle A_i,H\rangle \leq 0,~i \in \mathcal{T}_*;~\langle A_i, H\rangle \in \mathbb{R},~i \notin \mathcal{T}_*\}$. Hence $0= {\rm d}\varphi(X^*)(H) = \underset{i \in \mathcal{T}_*}{\sum}z_i^*\langle  A_i,H\rangle.$ Since $\langle  A_i,H\rangle\leq 0$ and $ z_i^*\ge 0$ for $i \in \mathcal{T}_*$, we have $z_i^*\langle  A_i,H\rangle  = 0$,  $i \in \mathcal{T}_*$. It follows that  $H \in S$. Thus, $K_{\varphi}(X^*,V) \subseteq S$.
		On the other hand, for any $H \in S$, it is clear that $z_i^*\langle  A_i,H\rangle= 0$ for $i\in \mathcal{T}_*$. For $i \in \bar{\mathcal{T}}_*$, $z_i^* =0$ and $\langle  A_i,H\rangle \in \R$, we also have $z_i^*\langle  A_i,H\rangle =0$, $ i \in \bar{\mathcal{T}}_*$.  Therefore,
		\begin{equation*}
			\langle V,H \rangle =\underset{i \in \mathcal{T}_*}{\sum}z_i^*\langle  A_i,H\rangle +\underset{i \notin \mathcal{T}_*}{\sum}z_i^*\langle  A_i,H\rangle =0.
		\end{equation*}
		Since $H \in S$, it is obvious that $H \in {\rm dom~d}\varphi$. It follows from Proposition \ref{limsubdiff} that $0={\rm d}\varphi(X^*)(H)=\langle V,H \rangle$, which implies that $H \in K_{\varphi}(X^*,V)$. Hence, $S \subseteq K_{\varphi}(X^*,V)$. In conclusion, $S = K_{\varphi}(X^*,V)$.
		
		\section*{Proof of  Proposition \ref{Prop 3.6}}
		It follows from \cite[Proposition 2]{L} , we have that for any $h \in {\rm dom}~{\rm d}g(u^* )$
		\begin{equation}\label{parabolic}
			{\rm d^2}g(u^*)(h,w) =0, ~ \forall w \in \{w \in \mathbb{R}^n~|~w_i \le 0 ~{\rm if} ~u_i^* = 0~{\rm and}~ h_i =0 \}.
		\end{equation}
		It follows from  \eqref{parabolic} and \cite[Exercise 13.63]{RW}  that for any $H \in {\rm dom}~{\rm d} \varphi(X^* )$, one has
		$$	{\rm d^2}\varphi(X^*)(H,W) = 	{\rm d^2}g(F(X^*))(AH,AW) =0,$$
		for $ W$ satisfying $\langle A_i,W \rangle\leq 0, \langle A_i,H \rangle =  0, ~i \in \mathcal{T}_*;~\langle A_i,W \rangle\in \mathbb{R},~i \notin \mathcal{T}_*$ and that $\varphi$ is parabolically epi-differentiable at $X^*$ for $H$ since the function $g$ is subdifferentially regular at $F(X^*)$ and parabolically epi-differentiable at $F(X^*)$ for $AH$.
		
		\section*{Proof of  Proposition \ref{second-derivative}}
		It follows from  \cite[Theorem 13.14]{RW} that 
		\begin{equation*}
			{\rm d}^2\varphi(X^*,V)(H)\geq \sup_{z \in \varLambda(X^*,V)}{\rm d}^2g(F(X^*),z)(AH)
		\end{equation*}
		always holds. It is enough to show that
		\begin{equation*}
			{\rm d}^2\varphi(X^*,V)(H)\leq \sup_{z \in \varLambda(X^*,V)}{\rm d}^2g(F(X^*),z)(AH).
		\end{equation*}
		For any $H \in K_{\varphi}(X^*,V)$, we have $\mbox{d}\varphi(X^*)(H)=\langle V,H\rangle$. Let ${z} \in \varLambda(X^*,V)$ be such that ${z} \in \partial g(F(X^*)) =\partial^p g(F(X^*))$ and $A^T{z}=V$. Then 
		\begin{equation*}
			{\rm d} g(F(X^*))(AH) =\mbox{d}\varphi(X^*)(H) =\langle V,H\rangle=\langle {z},AH\rangle,
		\end{equation*}
		i.e., $AH \in K_g(F(X^*),z)$. By Lemma \ref{Prop 3.4}, there exist some $\bar{W}$ such that
		\begin{align}\label{eq3.2}
			\mbox{d}^2 g(F(X^*),z)(AH)
			&= {\rm d^2} g(F(X^*))(AH,A\bar{W})-\langle z,A\bar{W}\rangle \notag\\
			&={\rm d^2}\varphi(X^*)(H,\bar{W})-\langle V,\bar W\rangle.
		\end{align}
		Let $\tilde{F}: \mathbb{R}^{m \times p}\times\mathbb{R} \to \mathbb{R}^n\times\mathbb{R}$, $\tilde{F}(X,r)=(F(X),r)$.  $\tilde{F}$ is an affine mapping since $F(X)$ is an affine mapping. We have
		\begin{equation*}
			\mbox{epi}\varphi=\tilde{F}^{-1}\mbox{epi}g.
		\end{equation*}
		As a union of finitely many convex  polyhedral sets $\mbox{epi}\varphi$ is second order regular at $X^*$ in the direction $H \in K_\varphi(X^*,V)$ for $V\in \partial \varphi(X^*)$  by \cite[Page 203 ]{BS}.  Hence $\varphi$ is second order regular by \cite[Definition 3.93]{BS}. It follows from \cite[Proposition 3.4]{Twice} that $\mbox{d}^2 \varphi(X^*,V)(H) >-\infty$. Thus, $\varphi$ is parabolically regular at $X^*$ for $V = \partial^p\varphi(X^*) =\partial\varphi(X^*)$ in the direction $H$ by \cite[Proposition 3.103]{BS}. Thus, 
		\begin{align*}
			\mbox{d}^2 \varphi(X^*,V)(H)&=\inf_{W} \{{\rm d^2}\varphi(X^*)(H,W)-\langle V,W\rangle\}\\
			& \leq {\rm d^2}\varphi(X^*)(H,\bar W)-\langle V,\bar{W}\rangle \\
			& =\mbox{d}^2 g(F(X^*),z)(AH)\\
			& \leq \sup_{z\in \varLambda(X^*,V)} \mbox{d}^2g(F(X^*),z)(AH),
		\end{align*}
		where the second equality used \eqref{eq3.2}. Therefore, for any~$H \in K_{\varphi}(X^*,V)$, we have 
		\begin{equation*}
			\mbox{d}^2\varphi(X^*,V)(H)=\sup_{z\in \varLambda(X^*,V)}\mbox{d}^2g(F(X^*),z)(AH),
		\end{equation*}
		and the supremum can be attained at some $\hat{z} \in \varLambda(X^*,V)$. Combining with \cite[Theorem 3] {L} and Proposition \ref{Prop 3.4}, we have 
		\begin{align*} 
			&{\rm d}^2\varphi(X^*,V)(H)\\
			=& \mbox{d}^2g(F(X^*),\hat{z})(AH)\\
			= &\left\{
			\begin{array}{clc}
				0,   &~\hat{z}_i\langle  A_i,H \rangle=0,~
				\langle  A_i,H \rangle\leq 0,~\hat{z}_i\geq 0,~i\in \mathcal{T}_*,\\
				+\infty, &~\mbox{ otherwise.} 
			\end{array}
			\right.
		\end{align*}
		
		\section*{Proof of  Theorem \ref{Second-condition}}
		It follows from \cite[Excercise 13.18]{RW} and Proposition \ref{second-derivative} that 
		\begin{align*}
			&\mbox{d}^2(f+\lambda \varphi)(X^*, 0)(D) =\langle D, \nabla^2f(X^*)D \rangle+ {\rm d^2} \varphi(X^*, -\nabla f(X^*))(D),
		\end{align*}
		for any $D \in \left\{D  \in \R^{m\times p}~|~z^*_i\langle D,A_i \rangle =0,~
		\langle D,A_i \rangle \leq 0,~z^*_i\geq 0,~~i\in \mathcal{T}_*\right\}$. Using  \eqref {critical-varphi}, we obtain 
		\begin{align*}
			&{\rm d^2} \varphi(X^*, -\nabla f(X^*))(D)\\
			=&\left\{
			\begin{array}{clc}
				0,  ~& z^*_i\langle D,A_i \rangle =0,~
				\langle D,A_i \rangle \leq 0,~z^*_i\geq 0,~~i\in \mathcal{T}_*,\\
				+\infty, ~&  \mbox{ otherwise }.
			\end{array}
			\right.
		\end{align*}
		Hence ,
		\begin{align*}
			&\mbox{d}^2(f+\lambda \varphi)(X^*, 0)(D)\\
			=&\left\{
			\begin{array}{clc}
				\langle D, \nabla^2f(X^*)D \rangle,  ~& D \in \left\{\begin{aligned}
					&D  \in \R^{m\times p}~|~z^*_i \langle D,A_i \rangle =0,\\
					&\langle D,A_i \rangle  \leq 0,~z^*_i\geq 0,~~i\in \mathcal{T}_*
				\end{aligned}\right\},\\
				+\infty, ~&  \mbox{ otherwise }.
			\end{array}
			\right.
		\end{align*}
		The conclusions then follow from\cite[Theorem 13.24]{RW}.  		
		
		\section*{Proof of  Theorem \ref{stationary-equation}}
		$(\romannumeral 1)$ It follows from Fermat's rule \cite[Theorem 10.1]{RW}	that $-\nabla f(X^*) \in \lambda \partial \varphi(X^*)$ if $X^*$ is a local minimizer of \eqref{Prim}. By using the fact $\lambda \partial \|(\cdot)_+\|_0 = \partial \|(\cdot)_+\|_0$ and Proposition \ref{limsubdiff}, we have that 
		$$-\nabla f(X^*) \in A^T\partial\|F(X^*)_+\|_0.$$ 
		Thus, there exists $z^* \in \partial \|F(X^*)_+\|_0 $ such that 
		\begin{equation}\label{eq2.5}	
			\nabla f(X^*) +A^Tz^*=0,
		\end{equation}
		where $z_i^*\geq 0 $ for $ i \in \mathcal{T}_*$ and $z^*_i=0 $ for $ i\notin  \mathcal{T}_*.$ 	
		From $0<\alpha <\alpha_*$, we obtain
		$$F_i(X^*) \geq \min_{i:F_i(X^*)>0} F_i(X^*) = \sqrt{2\lambda\alpha_1}\geq \sqrt{2\lambda\alpha_*}>\sqrt{2\lambda\alpha}$$
		and 
		$$z_i^*\leq \min_{i:z_i^*>0} z_i^* = \sqrt{2\lambda\alpha_2}\leq \sqrt{2\lambda\alpha_*}<\sqrt{2\lambda\alpha}.$$
		Therefore, we have that 
		\begin{equation*}\label{g_i(X^*)}
			\left\{\begin{aligned}
				&F_i(X^*) \in (-\infty, 0)\cup(\sqrt{2\lambda\alpha}, +\infty), & z^*_i = 0,\\
				&0<z_i^*<\sqrt{2\lambda\alpha} ~{\rm or}~ z^*_i =0.,~&F_i(X^*) = 0.
			\end{aligned}
			\right.
		\end{equation*}
		It follows that  $F(X^*) \in {\rm prox}_{\alpha \lambda \|(\cdot)_+\|_0}(F(X^*) + \alpha z^*)$.
		This together with \eqref{eq2.5}, implies that $(X^*, z^*)$ is a P-stationary point for problem \eqref{Prim}.\\
		$(\romannumeral 2)$A P-stationary point $(X^*,z^*)$  satisfies \eqref{Proximal2}, which implies that $z^* \in \partial \|F(X^*)_+\|_0 $. Consequently,  we can obtain
		$$0 = \nabla f(X^*) + A^Tz^* \in  \nabla f(X^*) + A^T\partial\|(F(X^*)_+\|_0,$$ i.e., $X^*$ is a local minimizer  of problem \eqref{Prim}.
		
		\section*{Proof of  Lemma \ref{4.1}}
		Since $(X^*,z^*)$ is a P-stationary point of  \eqref{Prim}, we have  $\Psi(X^*,z^*;\mathcal{T}_*)=0$, i.e.,
		\begin{equation*}
			\nabla f(X^*) + A^T_{\mathcal{T}_*}z_{\mathcal{T}_*} =0,~u_{\mathcal{T}_*}^*=0 ~{\rm and }~z_{\bar{\mathcal{T}}_*}^*=0.
		\end{equation*}
		For $0<\alpha<\alpha_*$ and $\mathcal{T}_*=\mathcal{T}_o^*\cup \mathcal{T}_1^*$, we  note that $\mathcal{T}_2^* \setminus\mathcal{T}_o^* \subseteq \bar{\mathcal{T}}_*$ and then
		$$z^*_{\mathcal{T}_2^* \setminus\mathcal{T}_o^*}=0.$$
		From the proof of Theorem \ref{stationary-equation}, we obtain
		$$\left\{\begin{aligned}
			&u_i^* >\sqrt{2\alpha\lambda},~& u_i^*>0,\\
			&z_i^* < \sqrt{2\lambda/\alpha},~& z_i^*>0.
		\end{aligned}
		\right.$$ 
		Therefore, we have the following facts:
		\begin{align*}
			\mathcal{T}_o^* &= \{i~|~u_i^* = 0,~\alpha z_i^* \in \{0, \sqrt{2\alpha\lambda}\}\}\\
			& = \{i~|~u_i^* = 0,~z_i^* =0\},\\
			\mathcal{T}_2^*\setminus\mathcal{T}_o^* &= \{i~|~u_i^* \neq 0,~u_i^* + \alpha z_i^* =\sqrt{2\alpha \lambda}\}\\
			& = \{i~|~u_i^*=\sqrt{2\alpha \lambda}\}\\
			& = \emptyset,
		\end{align*}
		these facts lead to
		$$	\mathcal{T}_2^*  =	\mathcal{T}_o^* \cup (\mathcal{T}_2^*\setminus\mathcal{T}_o^*) = \mathcal{T}_o^*=\{i~|~u_i^* = 0,~z_i^* =0\}$$
		and
		$$\mathcal{T}_* = \mathcal{T}_1^*\cup \mathcal{T}_o^*= \mathcal{T}_1^*\cup \mathcal{T}_2^*,~ \bar{\mathcal{T}}_* = \mathcal{T}_3^*.$$
		We define a neighborhood $B((X^*,z^*), \delta_1^*)$ of $(X^*,z^*)$ with 
		$$0<\delta_1^*= \min\{ \hat{\delta}^*, \beta_1,~\beta_2\},$$	
		where $\beta_1 =\frac{1}{\|A\|_F+\alpha}\{\sqrt{\lambda\alpha/2}- \underset{j \in \mathcal{T}_1^*}{\max} |u_j^* + \alpha z_j^* -\sqrt{\lambda\alpha/2}|\}$ and $\beta_2 =\frac{1}{\|A\|_F+\alpha}\{\underset{j \in \mathcal{T}_3^*}{\min} |u_j^* + \alpha z_j^* -\sqrt{\lambda\alpha/2}|-\sqrt{\lambda\alpha/2}\}$. For any $(X,z) \in B((X^*,z^*), \delta_1^*)$  and $i \in \mathcal{T}_1^*$, we have that
		\begin{align*}
			&|u_i + \alpha z_i -\sqrt{\lambda\alpha/2}| \\
			\le &|u_i + \alpha z_i - (u_i^* + \alpha z_i^*)|+ |u_i^* + \alpha z_i^* -\sqrt{\lambda\alpha/2}|\\
			<&(\|A\|_F+\alpha)\delta_1^* + |u_i^* + \alpha z_i^* -\sqrt{\lambda\alpha/2}|\\
			\leq &\sqrt{\lambda\alpha/2}-\underset{j \in \mathcal{T}_1^*}{\max} |u_j^* + \alpha z_j^* -\sqrt{\lambda\alpha/2}| + \underset{i \in \mathcal{T}_1^*}{\max} |u_i^* + \alpha z_i^* -\sqrt{\lambda\alpha/2}|\\
			= &\sqrt{\lambda\alpha/2},
		\end{align*} 
		Therefore,  $u_i + \alpha z_i  \in (0,\sqrt{2\lambda\alpha} )$, i.e., $i\in \mathcal{T}_1$ and $\mathcal{T}_1^* \subseteq \mathcal{T}_1.$
		Similarly, for any $(X,z) \in B((X^*,z^*), \delta_1^*)$  and $i \in \mathcal{T}_3^*$, we have that	
		\begin{align*}
			&|u_i + \alpha z_i -\sqrt{\lambda\alpha/2}| \\
			\ge& |u_i^* + \alpha z_i^* -\sqrt{\lambda\alpha/2}|- |u_i + \alpha z_i - (u_i^* + \alpha z_i^*)| \\
			> &|u_i^* + \alpha z_i^* -\sqrt{\lambda\alpha/2}|- (\|A\|_F+\alpha)\delta_1^*  \\
			\ge& |u_i^* + \alpha z_i^* -\sqrt{\lambda\alpha/2}| - \underset{j \in \mathcal{T}_3^*}{\min} |u_j^* + \alpha z_j^* -\sqrt{\lambda\alpha/2}|+\sqrt{\lambda\alpha/2}\\
			\ge &\underset{i \in \mathcal{T}_3^*}{\min} |u_i^* + \alpha z_i^* -\sqrt{\lambda\alpha/2}|-\underset{j \in \mathcal{T}_3^*}{\min} |u_j^* + \alpha z_j^* -\sqrt{\lambda\alpha/2}| +\sqrt{\lambda\alpha/2}\\
			\ge& \sqrt{\alpha \lambda/2	}.
		\end{align*} 
		Hence,  we prove that $i\in \mathcal{T}_3$ and $\mathcal{T}_3^* \subseteq \mathcal{T}_3.$
		From the above conclusions, we also have 
		$$\mathcal{T}_2= \mathbb{N}_n\setminus (\mathcal{T}_1\cup \mathcal{T}_3)\subseteq \mathbb{N}_n\setminus (\mathcal{T}_1^*\cup \mathcal{T}_3^*) =  \mathcal{T}_2^*.$$
		In conclusion, for any $(X,z) \in B((X^*,z^*), \delta_1^*)$, it holds that
		$$\mathcal{T}_1^* \subseteq \mathcal{T}_1,~\mathcal{T}_3^* \subseteq \mathcal{T}_3,~\mathcal{T}_2 \subseteq \mathcal{T}_2^*.$$
		By these relations, we can deduce that for a sufficiently small  $\delta_1^*$, 
		$$\mathcal{T}_1\setminus \mathcal{T}_1^*\subseteq \mathcal{T}_2^*,~\mathcal{T}_3\setminus \mathcal{T}_3^*\subseteq \mathcal{T}_2^*.$$
		Now we can show that
		\begin{align*}
			\mathcal{T}=\mathcal{T}_o\cup \mathcal{T}_1 &=\mathcal{T}_o\cup \mathcal{T}_1^*\cup (\mathcal{T}_1\setminus \mathcal{T}_1^*)&\\
			&\subseteq\mathcal{T}_2^*\cup \mathcal{T}_1^*\cup \mathcal{T}_2^*~&({\rm by}~\mathcal{T}_o \subseteq \mathcal{T}_2 \subseteq \mathcal{T}_2^* )\\
			&= \mathcal{T}_1^*\cup \mathcal{T}_2^*.
		\end{align*}
		It follows from $(X^*, z^*)$ being a P-stationary point and $\mathcal{T}\subseteq \mathcal{T}_*$ that $\Psi(X^*, z^*; \mathcal{T})=0.$
		
		\section*{Proof of Lemma \ref{4.2}}
		For any given two matrices $D'$ and $D$, we have 
		\begin{align}\label{matrix}
			\sigma_{\rm max}(D'-D)&\ge \underset{i}{\rm max}|\sigma_i(D')-\sigma_i(D)|\ge |\sigma_{i_0}(D')-\sigma_{i_0}(D)|\notag\\
			& \ge \sigma_{i_0}(D')-\sigma_{\rm min}(D)\ge \sigma_{\rm min}(D')-\sigma_{\rm min}(D),
		\end{align}
		where the first inequality is from \cite[Page 76]{HM} and $i_0$ satisfies $\sigma_{i_0}(D)=\sigma_{\rm min}(D)$. From \cite[Page 80]{HM}, we also have
		\begin{equation}\label{max}
			\sigma_{\rm max}(D')\le \sigma_{\rm max}(D'-D) + \sigma_{\rm max}(D).
		\end{equation}
		Since  $\mathcal{T}\subseteq \mathcal{T}_*$,   $H(\mathcal{T})$ is a submatrix of $H(\mathcal{T}_*)$. Hence
		$$\sigma_{\rm max}(H(\mathcal{T}))\le \sigma_{\rm max}(H(\mathcal{T}_*)),~ \sigma_{\rm min}(H(\mathcal{T}))\ge \underset{\mathcal{T}\subseteq \mathcal{T}_*}{\min}\sigma_{\rm min}(H(\mathcal{T})).$$
		By \cite[Proposition 10]{L}, the second-order sufficient condition implies that  $\nabla \tilde{\Psi}(X^*,z^*;\mathcal{T}_*)$ is non-singular. Since for  any $(X,z)\in B((X^*,z^*),\delta^*_2)$ with $\delta^*_2 \le\delta^*_1$, we have $\mathcal{T} \subseteq \mathcal{T}_*$ by Lemma \ref{4.1} and hence
		$$\nabla \tilde{\Psi}(X^*,z^*;\mathcal{T}) = \begin{pmatrix}
			H(\mathcal{T})&0\\0&I
		\end{pmatrix}$$ 
		is also non-singular for any $\mathcal{T}\subseteq \mathcal{T}_*$. Thus, $H(\mathcal{T})$ is non-singular and $\sigma_{\rm min}(H(\mathcal{T}))>0$ for any $\mathcal{T}\subseteq \mathcal{T}_*$. 
		Next we have the fact
		\begin{align}\label{c_*}
			\sigma_{\rm min}(\nabla \tilde{\Psi}(X^*,z^*;\mathcal{T})) &= {\rm min}\{1,\sigma_{\rm min}(H(\mathcal{T}))\} \notag\\
			&\ge {\rm min}\{1, \underset{\mathcal{T}\subseteq \mathcal{T}_*}{\rm min}{\sigma_{\rm min}}(H(\mathcal{T}))\} =2c_*,
		\end{align}
		\begin{align}\label{C_*}
			\sigma_{\rm max}(\nabla \tilde{\Psi}(X^*,z^*;\mathcal{T})) &= {\rm max}\{1,\sigma_{\rm max}(H(\mathcal{T}))\}\notag \\
			&\le {\rm max}\{1,\sigma_{\rm max}(H(\mathcal{T}_*))\}=\frac{C_*}{2}.
		\end{align}
		The locally Lipschitz continuity of $\nabla^2f$ around $X^*$ with $L_*$ yields 
		\begin{align}\label{21}
			\|\nabla \tilde{\Psi}(X^*,z^*;\mathcal{T})-\nabla \tilde{\Psi}(X,z;\mathcal{T})\| \overset{\eqref{spectral-norm}}{\le} &\|\nabla^2 f(X^*)-\nabla^2 f(X)\|_F\notag\\
			\le& L_*\|X-X^*\|_F\notag \\
			\le & L_*\delta_2^*\notag.\\
			\overset{\eqref{parameter2}}{\le} & \frac{c_*}{2},
		\end{align}
		Now by \eqref{matrix}, \eqref{c_*} and \eqref{21}, we have that
		{\begin{align*}
				&\sigma_{\min}(\nabla\tilde{\Psi}_{\gamma}(X,z;\mathcal{T}))\\
				\overset{\eqref{matrix}}{\ge}& \sigma_{\min}(\nabla\tilde{\Psi}(X,z;\mathcal{T}))-\sigma_{\max}(\nabla\tilde{\Psi}(X,z;\mathcal{T})-\nabla\tilde{\Psi}_{\gamma}(X,z;\mathcal{T}))\\
				=&  \sigma_{\min}(\nabla \tilde{\Psi}(X,z;\mathcal{T}))-\gamma\\
				\overset{\eqref{matrix}}{\ge}& \sigma_{\min}(\nabla\tilde{\Psi}(X^*,z^*;\mathcal{T}))-\gamma\\
				~&-\sigma_{\max}(\nabla\tilde{\Psi}(X^*,z^*;\mathcal{T})-\nabla\tilde{\Psi}(X,z;\mathcal{T}))\\
				\overset{\eqref{21}}{\ge}& \sigma_{\min}(\nabla\tilde{\Psi}(X^*,z^*;\mathcal{T}))-\frac{c_*}{2}-\gamma\\
				\overset{\gamma \le \frac{c_*}{2}}{\ge}& \sigma_{\min}(\nabla\tilde{\Psi}(X^*,z^*;\mathcal{T}))-c_*\\
				\overset{\eqref{c_*}}{\ge}&  c_*.
		\end{align*}}
		It follows from \eqref{max}, \eqref{C_*} and \eqref{21}, we also have 
		\begin{align*}
			~&\sigma_{\max}(\nabla\tilde{\Psi}_{\gamma}(X,z;\mathcal{T}))\\
			\overset{\eqref{max}}{\le}& \sigma_{\max}(\nabla\tilde{\Psi}(X,z;\mathcal{T})-\nabla\tilde{\Psi}_{\gamma}(X,z;\mathcal{T}))+ \sigma_{\max}(\nabla\tilde{\Psi}(X,z;\mathcal{T}))\\
			=& \gamma + \sigma_{\rm max}(\nabla\tilde{\Psi}(X,z;\mathcal{T}))\\
			\overset{\eqref{max}}{\le}& \sigma_{\rm max}(\nabla\tilde{\Psi}(X,z;\mathcal{T})-\nabla\tilde{\Psi}(X^*,z^*;\mathcal{T}))+\gamma \\
			~&+ \sigma_{\rm max}(\nabla\tilde{\Psi}(X^*,z^*;\mathcal{T}))\\
			\overset{\eqref{21}}{\le}&  \sigma_{\max}(\nabla\tilde{\Psi}(X^*,z^*;\mathcal{T}))+\frac{c_*}{2}+\gamma \\
			\overset{\eqref{C_*}}{\le}&  \frac{C_*}{2} +c_*\\
			{\le}& C_*.
		\end{align*}
		The last inequality uses \eqref{parameter1}, i.e., $c_* \le \frac{C_*}{2}$. This complete the proof.
		
		\section*{Proof of Theorem \ref{convergence}}	
		It follows from Lemma \ref{4.1} and the fact $\delta^*\le \delta_1^*\le \delta^*_2$  that for the index set $\mathcal{T}_0$ corresponding to  $(X^0,z^0) \in B((X^*,z^*),\delta^*),$ we have 
		$$\tilde{\Psi}(X^*,z^*;\mathcal{T}_0) = {\rm vec}(\Psi(X^*,z^*;\mathcal{T}_0)) =0.$$
		By using the fact $\gamma_0\leq \gamma_{-1}\leq \frac{c_*}{2}$ and the Newton equation, we have 
		$$\nabla\tilde{\Psi}_{\gamma_0}(X^0,z^0;\mathcal{T}_0) {\rm vec} (D^0)=-\tilde{\Psi} (X^0,z^0;\mathcal{T}_0) .$$
		According to  the nonsingularity of $\nabla\tilde{\Psi}_{\gamma_0}(X^0,z^0;\mathcal{T}_0)$, it is clear that ${\rm vec} (D^0)$ is well defined.
		For convenience, let $(\tilde{x},z) = ({\rm vec}(X),z)$.
		We have the following conclusion
		\begin{align*}
			&\|(\tilde{x}^1,z^1) -(\tilde{x}^*,z^*)\| \\
			=&\|(\tilde{x}^0,z^0) -(\tilde{x}^*,z^*)-\nabla \tilde{\Psi}_{\gamma_0}(X^0,z^0;\mathcal{T}_0)^{-1}\tilde{\Psi} (X^0,z^0;\mathcal{T}_0)\|\\
			\le &\|\nabla \tilde{\Psi}_{\gamma_0}(X^0,z^0;\mathcal{T}_0)\|^{-1}\|\nabla \tilde{\Psi}_{\gamma_0}(X^0,z^0;\mathcal{T}_0)[(\tilde{x}^0,z^0)-(\tilde{x}^*,z^*)]\\
			~&-\tilde{\Psi} (X^0,z^0;\mathcal{T}_0) + \tilde{\Psi}(X^*,z^*;\mathcal{T}_0)\|\\
			\overset{\eqref{15}}{\le}& \frac{1}{c_*}\int_{0}^{1}\|\nabla \tilde{\Psi}_{\gamma_0}(X^0,z^0;\mathcal{T}_0)-\nabla \tilde{\Psi}((X^*,z^*)\\
			&+s((X^0,z^0)-(X^*,z^*));\mathcal{T}_0)\| \|(\tilde{x}^0,z^0) -(\tilde{x}^*,z^*)\|{\rm d}s
		\end{align*}
		and
		\begin{eqnarray*}
			&&\|\nabla \tilde{\Psi}_{\gamma_0}(X^0,z^0;\mathcal{T}_0)-\nabla \tilde{\Psi}((X^*,z^*)+s((X^0,z^0) -(X^*,z^*)));\mathcal{T}_0)\|\\
			&\overset{\eqref{spectral-norm}}{\le}&\|\nabla^2 f(X^0) - \nabla^2 f(X^*+s(X^0-X^*))\|_F+\gamma_0\\
			&\overset{\eqref{gamma}}{\le}& L_*(1-s)\|X^0-X^*\|_F +\rho \|\tilde{\Psi}(X^0,z^0;\mathcal{T}_0)\|\\
			&\overset{\eqref{FnormSubtraction}}{=}& L_*(1-s)\|\tilde{x}^0-\tilde{x}^*\| +\rho \|\tilde{\Psi}(X^0,z^0;\mathcal{T}_0)\|\\
			&\le& L_*(1-s)\|(\tilde{x}^0,z^0)-(\tilde{x}^*,z^*)\| +\rho \|\nabla \tilde{\Psi}_{\gamma_0}(X^0,z^0;\mathcal{T}_0){\rm vec} (D^0)\|\\
			&\overset{\eqref{15}}{\le} & L_*(1-s)\|(\tilde{x}^0,z^0)-(\tilde{x}^*,z^*)\| +\rho C_*\|(\tilde{x}^1,z^1)-(\tilde{x}^0,z^0)\|\\
			&\le& (L_*(1-s)+\rho C_*)\|(\tilde{x}^0,z^0)-(\tilde{x}^*,z^*)\|+\rho C_* \|(\tilde{x}^1,z^1)-(\tilde{x}^*,z^*)\|.
		\end{eqnarray*}
		By the above inequalities, we obtain 
		\begin{eqnarray*}
			&&\|(\tilde{x}^1,z^1)-(\tilde{x}^*,z^*)\|\\
			&\le& \frac{1}{c_*}\|\int_{0}^{1}(L_*(1-s)+\rho C_*)\|(\tilde{x}^0,z^0)-(\tilde{x}^*,z^*)\|^2\\
			&~&+\rho C_* \|(\tilde{x}^1,z^1)-(\tilde{x}^*,z^*)\|\|(\tilde{x}^0,z^0)
			-(\tilde{x}^*,z^*)\|{\rm d}s\\
			&=& \frac{L_*+2\rho C_*}{2c_*}\|(\tilde{x}^0,z^0)-(\tilde{x}^*,z^*)\|^2+\frac{\rho C_*}{c_*} \|(\tilde{x}^1,z^1)-(\tilde{x}^*,z^*)\|\\
			&~&\|(\tilde{x}^0,z^0)-(\tilde{x}^*,z^*)\|\\
			&\le& \frac{\theta_*}{2}\|(\tilde{x}^0,z^0)-(\tilde{x}^*,z^*)\|^2+\frac{\rho C_*\delta^*}{c_*} \|(\tilde{x}^1,z^1)-(\tilde{x}^*,z^*)\|\\
			&\overset{\eqref{delta2}}{\le}&\frac{\theta_*}{2}\|(\tilde{x}^0,z^0)-(\tilde{x}^*,z^*)\|^2 +\frac{1}{2}\|(\tilde{x}^1,z^1)-(\tilde{x}^*,z^*)\|,
		\end{eqnarray*}
		where $\theta _*=\frac{L_*+2\rho C_*}{c_*}$. We have $$\|(\tilde{x}^1,z^1)-(\tilde{x}^*,z^*)\|\le {\theta_*}\|(\tilde{x}^0,z^0)-(\tilde{x}^*,z^*)\|^2.$$
		It follows from \eqref{FnormSubtraction} that the above inequality is equivalent to
		\begin{equation}
			\|X^1-X^*\|_F+\|z^1-z^*\|\le \frac{L_*+2\rho C_*}{c_*}(\|X^0-X^*\|_F^2+\|z^0-z^*\|^2).
		\end{equation}
		It follows from $\delta^*\le \frac{1}{2\theta_*}$ and $(X^0,z^0) \in B((X^*,z^*),\delta^*)$ that
		\begin{align*}
			\|X^1-X^*\|_F+\|z^1-z^*\|&\le \theta_*(\|X^0-X^*\|_F^2+\|z^0-z^*\|^2)\\
			&\le \delta^*\theta_*(\|X^0-X^*\|_F+\|z^0-z^*\|)\\
			&\le \frac{1}{2}(\|X^0-X^*\|_F+\|z^0-z^*\|)\\
			&< \delta^*,	
		\end{align*}
		which means $(X^1,z^1)\in B((X^*,z^*),\delta^*)$. Replacing $\mathcal{T}_0$ by $\mathcal{T}_1$ and noting $\gamma_1\leq\gamma_0 \leq \frac{c_*}{2}$, we can deduce similarly that
		\begin{align*}
			\|X^2-X^*\|_F+\|z^2-z^*\|&\le\theta_*(\|X^1-X^*\|_F^2+\|z^1-z^*\|^2) \\
			&\le\frac{1}{2}(\|X^1-X^*\|_F+\|z^1-z^*\|),
		\end{align*}
		${\rm vec}(D^1)$ is also well defined and $(X^2,z^2) \in B((X^*,z^*),\delta^*)$. By the induction, we can conclude that ${\rm vec}(D^k)$ is well defined and $(X^{k+1},z^{k+1}) \in B((X^*,z^*),\delta^*)$ be such that
		\begin{align*}
			\|X^{k+1}-X^*\|_F + \|z^{k+1}-z^*\| \le&\frac{L_*+2\rho C_*}{c_*}(\|X^k-X^*\|_F^2 \\
			&+\|z^k-z^*\|^2).
		\end{align*}
		Thus, we have $(X^k,z^k)\to (X^*,z^*)$ and
		\begin{align*}
			{\rm vec}(D^{k+1}) &=(\tilde{x}^{k+1},z^{k+1})-(\tilde{x}^k,z^k)\\ 
			&=(\tilde{x}^{k+1},z^{k+1})-(\tilde{x}^*,z^*)+(\tilde{x}^*,z^*)-(\tilde{x}^k,z^k) \to 0,
		\end{align*}
		which implies $\underset{k \to \infty}{\lim } (D^k_{X}, D^k_{\mathcal{T}_k}, D^k_{{\bar{\mathcal{T}}_k}}) =0.$
		This completes the proof of  $(\romannumeral 1)$ and $(\romannumeral 2)$.\\
		$(\romannumeral 3)$ The above proof shows that $(X^k,z^k) \in B((X^*,z^*),\delta^*)$. It follows from Lemma  \ref{4.1} that 
		$$ \tilde{\Psi}(X^*,z^*;\mathcal{T}_k) =0.$$
		Let $(X^k,z^k)_\beta = (X^*,z^*)+\beta[(X^k,z^k)-(X^*,z^*)]$, where $\beta \in [0,1]$, we have $(X^k,z^k)_\beta  \in B((X^*,z^*),\delta^*) \subseteq  B((X^*,z^*),\delta^*_2)$. Set $\gamma=0$ and $(X,z) =(X^k,z^k)_\beta$ in \eqref{15},  we have that $\nabla \tilde{\Psi}((X^k,z^k)_\beta;\mathcal{T}_k)$ has the same bounds as in \eqref{15}, namely
		\begin{equation}
			C_*\ge \sigma_{\rm max}(\nabla \tilde{\Psi}((X^k,z^k)_\beta;\mathcal{T}_k))\ge \sigma_{\rm min}(\nabla \tilde{\Psi}((X^k,z^k)_\beta;\mathcal{T}_k) \ge c_*.
		\end{equation}
		It follows from  the differentiablity of $\tilde{\Psi}(\cdot,\mathcal{T}_k)$ for the fixed $\mathcal{T}_k$ and the Mean-value theorem that there exists a $\beta_0\in (0,1)$ satisfying
		\begin{align*}
			&\|\Psi(X^k,z^k;\mathcal{T}_k)\|_F\\
			=& \|\tilde{\Psi}(X^k,z^k;\mathcal{T}_k)\| \\
			=& \|\tilde{\Psi}(X^*,z^*;\mathcal{T}_k)+\nabla \tilde{\Psi} ((X^k,z^k)_{\beta_0};\mathcal{T}_k)[(\tilde{x}^k,z^k)-(\tilde{x}^*,z^*)] \|\\
			=& \|\nabla \tilde{\Psi}((X^k,z^k)_{\beta_0};\mathcal{T}_k)[(\tilde{x}^k,z^k)-(\tilde{x}^*,z^*)]\|\\
			\in& [c_*\|(\tilde{x}^k,z^k)-(\tilde{x}^*,z^*)\|, C_* \|(\tilde{x}^k,z^k)-(\tilde{x}^*,z^*)\|].
		\end{align*}
		The third equality uses the fact that  $\tilde{\Psi}(X^*,z^*;\mathcal{T}_k)=0$. Therefore,
		\begin{align*}
			\|\Psi(X^k,z^k;\mathcal{T}_k)\|_F &\le C_* \|(\tilde{x}^k,z^k)-(\tilde{x}^*,z^*)\|\\
			&\le \theta_*C_*\|(\tilde{x}^{k-1},z^{k-1})-(\tilde{x}^*,z^*)\|^2\\
			& \le \theta_* C_* /c_*^2\|\Psi(X^{k-1},z^{k-1};\mathcal{T}_{k-1})\|_F^2\\
			&\le \theta_* C_*^3 /c_*^2\|(\tilde{x}^{k-1},z^{k-1})-(\tilde{x}^*,z^*)\|^2\\
			&\le \theta_* C_*^3 /c_*^2)2^2\|(\tilde{x}^{k-2},z^{k-2})-(\tilde{x}^*,z^*)\|^2\\
			&~\vdots\\
			&\le \theta_* C_*^3 /c_*^22^{2k-2}\|(\tilde{x}^0,z^0)-(\tilde{x}^*,z^*)\|^2.	
		\end{align*}
		We can verify that $\|\Psi((X^k,z^k);\mathcal{T}_k)\|_F<\epsilon$ if $k$ satisfies 
		\begin{align*}
			k\ge& {\rm log}_2(2\sqrt{\theta_*C_*^3/c_*^2}(\|X^0 -X^*\|_F+\|z^0-z^*\|))-{\rm log}_2\sqrt{\epsilon}\\
			= &{\rm log}_2(2\sqrt{\frac{(L_*+2\rho C_*)C_*^3}{c_*^3}}(\|X^0 -X^*\|_F+\|z^0-z^*\|))\\
			&-{\rm log}_2\sqrt{\epsilon}.
		\end{align*}
		This completes the whole proof.
	}
	

\end{document}